\documentclass{article}

\usepackage{amsmath,amssymb,amscd,amsthm}
\usepackage{diagrams}
\input cyracc.def

\newcommand{\QQ}{\mathbb Q}
\newcommand{\ZZ}{\mathbb Z}
\newcommand{\NN}{\mathbb N}

\newcommand{\Tau}{\frak{T}}
\renewcommand{\AA}{\mathbb A}
\newcommand{\Qp}{\QQ_p}
\newcommand{\Zp}{\ZZ_p}
\newcommand{\mup}{\mu_{p^\infty}}
\newcommand{\fK}{{\cal K}}
\newcommand{\fA}{{\frak A}}
\newcommand{\fV}{{\cal V}}
\newcommand{\fW}{{\cal W}}
\newcommand{\fH}{{\cal H}}
\newcommand{\fS}{{\cal S}}
\newcommand{\fM}{{\cal M}}
\newcommand{\fT}{{\cal T}}
\newcommand{\fU}{{\cal U}}
\newcommand{\fN}{{\cal N}}
\newcommand{\fX}{{\cal X}}
\newcommand{\calf}{\frak{f}}
\newcommand{\EE}{\mathbb E}
\newcommand{\FF}{\mathbb F}
\newcommand{\CC}{\mathbb C}
\newcommand{\pp}{{\frak p}}
\newcommand{\mm}{{\frak m}}

\newcommand{\Fp}{\FF_p}

\newcommand{\OfH}{O_{\fH}}
\newcommand{\EfK}{\EE_{\fK}}
\newcommand{\EfH}{\EE_{\fH}}
\newcommand{\EK}{\EE_{K}}
\newcommand{\EH}{\EE_{H}}
\newcommand{\fR}{{\cal R}}
\newcommand{\ffR}{{\cal R}_{\Qp}}
\newcommand{\fOK}{O_{\fK}}
\newcommand{\fOH}{O_{\fH}}
\newcommand{\muk}{\mu_{p^k}}
\newcommand{\tK}{\tilde{K}}
\newcommand{\vs}{\vspace{1ex}}
\newcommand{\bpi}{\bar{\pi}}

\newcommand{\fRKdaggern}{\fR_K^{\dagger,N}}
\newcommand{\fSKdaggern}{\fS_K^{\dagger,N}}

\newcommand{\fTKdaggern}{\fT_K^{\dagger,N}}
\newcommand{\fRdagger}{\fR^\dagger_{\Qp}}
\newcommand{\fRKdagger}{\fR_K^\dagger}
\newcommand{\whK}{\widehat{K}}
\newcommand{\OmegaRK}{\Omega^{n-1}_{\fR_K}}
\newcommand{\OmegaRKn}{\Omega^{n}_{\fR_K^+}}
\newcommand{\OmegaRKnlog}{\Omega^{n}_{\fR_K^+}(\log)}
\newcommand{\OmegaRKnloglog}{\Omega^{n}_{\fR_K^+}(\log)(\log)}
\newcommand{\del}{\frak{d}}

\newcommand{\OmegatwofRKlog}{\Omega^n_{\fR_K^+}(\log)}

\newcommand{\fRdaggern}{\fR^{\dagger,N}_{\Qp}}

\newcommand{\OmegatwofSKn}{\Omega^n_{\fSKdaggern}}
\newcommand{\OmegarfSKn}{\Omega^r_{\fSKdaggern}}
\newcommand{\OmegafSKnlog}{\Omega^1_{\fSKdaggern}(\log)}
\newcommand{\OmegatwofSKnlog}{\Omega^n_{\fSKdaggern}(\log)}
\newcommand{\OmegarfSKnlog}{\Omega^r_{\fSKdaggern}(\log)}

\newarrow{Into} C--->

\DeclareMathOperator{\coker}{coker}
\DeclareMathOperator{\trace}{trace}
\DeclareMathOperator{\Gal}{Gal}

\DeclareMathOperator{\image}{Im}
\DeclareMathOperator{\id}{id}

\DeclareMathOperator{\Tr}{Tr}
\newcommand{\ppp}{\frak{pr}}

\newtheorem{thm}{Theorem}
\newtheorem{prop}[thm]{Proposition}
\newtheorem{lem}[thm]{Lemma}
\newtheorem{cor}[thm]{Corollary}

\begin{document}

\title{Generalized logarithmic derivatives for $K_n$}

\author{Sarah Livia Zerbes}

\maketitle

\begin{center}
 {\it dedicated to John Coates}
\end{center}

\begin{abstract}
 \noindent We construct (generalized) logarithmic derivatives for general $n$-dimensional local fields 
 $\fK$ of mixed characteristics $(0,p)$ in which $p$ is not necessarily 
 a prime element with residue field $k$ such
 that $[k:k^p]=p^{n-1}$. For the construction of the logarithmic derivative map, we define $n$-dimensional 
 rings of 
 overconvergent series and show that - as in the $1$-dimensional case - they can be interpreted as 
 functions converging on some annulus of the open unit $p$-adic disc. Using the generalized logarithmic
 derivative, we give a new construction of Kato's $n$-dimensional dual exponential map.
 \vs


\end{abstract}

\tableofcontents

%++++++++++++++++++++++++++++++++++++++++++++++++++++++++++++++

 \section{Introduction}
 
  \subsection{The main results}
 
  Let $p$ be any prime number and $F$ a finite unramified extension of $\Qp$. For each $n\geq 1$, let
  $\mu_{p^n}$ be a primitive $p^nth$ root of unity such that $\mu_{p^n}^p=\mu_{p^{n-1}}$, and let
  $F_n=F(\mu_{p^n})$. Let $F_\infty=F(\mup)$ and $\Gamma=\Gal(F^{\text{cyc}}/F)$, where
  $F^{\text{cyc}}\subset F_\infty$
  is the cyclotomic $\Zp$-extension of $F$. One of the major tools in classical Iwasawa theory are 
  so-called Coleman power series (which were discovered by Coates and Wiles~\cite{coateswiles} and 
  Coleman~\cite{coleman}), which associate to every norm-compatible system of units $u=(u_n)$ in the tower 
  $F_\infty=\bigcup_nF_n$ a norm-invariant element $G_u(T)$ in the Iwasawa algebra
  $\Lambda(\Gamma)\cong O_F[[T]]$ such that 
  $G_u(\mu_{p^n})=u_n$. It was unknown for a long time how to generalize this construction for finite {\it
  ramified} extensions of $\Qp$, and only recently Cherbonnier and Colmez~\cite{cherbonniercolmez2}
  discovered a way to construct Coleman series for such extensions using a formal analogy between 
  $\Lambda(\Gamma)$ and Fontaine's rings of overconvergent series. 
  \vs
  
  In~\cite{fukaya2}, Fukaya extends the classical construction of Coleman power series to a certain class
  of $2$-dimensional local fields of mixed characteristics $(0,p)$ with imperfect residue fields. 
  Given such a field $\fH$, she replaces the cyclotomic $\Zp$-extension by a certain 
  $2$-dimensional false Tate curve extension $\fH_\infty$ of $\fH$ (the field 
  $\fH_\infty$ is the direct limit of a
  sequence of finite Galois extensions $\fH_n$ of $\fH$), and in analogy to the classical case,
  her Coleman isomorphism associates to every element in $\varprojlim K_2(O_{\fH_n})$ a norm-invariant 
  element in $K_2(O_{\fH}[[T]])$. Composing this map with the $d\log$-operator, she obtains a canonical
  homomorphism
  \begin{equation}
   \tau:\varprojlim K_2(O_{\fH_n})\rTo (\Omega^2_{O_{\fH}[[T]]})^{\psi=1},
  \end{equation}
  where $\psi$ is the derivative of the norm map. Because of its analogy with the classical
  case, this map is called the $2$-dimensional logarithmic derivative. 
  \vs
  
  However, Fukaya needs to assume for her construction that $p$ is a prime element in $\fH$. Using the
  theory of higher-dimensional fields of norms developed by Scholl in~\cite{scholl}, in this paper 
  we define
  $n$-dimensional rings of overconvergent series, and following the approach of Cherbonnier and Colmez
  in~\cite{cherbonniercolmez1}, we use these rings to extend Fukaya's construction of the logarithmic
  derivative map to the case of a ramified basefield
  and to $n$-dimensional local fields of mixed characteristic $(0,p)$.
  The main result of this paper is the following theorem:- 
  \vs
  
  \noindent {\bf Theorem~\ref{logderivrevisited}.}
  The logarithmic derivative gives a canonical homomorphism
  \begin{equation}
   \Tau:\varprojlim\whK_n(\fK_i)\rTo(\OmegatwofSKnlog)^{\psi=1}.
  \end{equation}

  In Section~\ref{relationdualexp}, we
  then use the logarithmic derivative map to give a new construction of Kato's dual exponential map for
  $K_n$.
  
%++++++++++++++++  
  
  \subsection{Acknowledgements}
  
   I am very grateful to John Coates for his interest and to him and Yann Weil for their 
   encouragement whenever I got fed up
   with the stubborn $K$-groups. Also, I would like to thank Max Karoubi for helpful discussions, Takako
   Fukaya for sending me copies of her paper~\cite{fukaya1} and Tony Scholl for giving me his 
   preprint~\cite{scholl} on higher-dimensional local fields.

%++++++++++++++++++++++++++++++++++++++++++++++++++++++++++++++

 \section{Setup}

%++++++++++++++++++

  \subsection{Fields of norms}
  
   Throughout this paper, we let $p$ be some prime $>2$. (The calculations in the
   case $p=2$ are slightly more difficult and will be dealt with in another paper.)
   \vs
   
   Let $K$ be a finite totally ramified extension of $\Qp$ with ring of integers $O_K$ (we allow the case
   when $K=\Qp$), and assume that
   $K\cap\QQ_p(\mu_{p^\infty})=\QQ_p$. We construct our $n$-dimensional local field inductively:- Let
   \begin{align}
    B_1 &= O_K[[T_1]][T_1^{-1}]^{\vee},\\
    B_i &= B_{i-1}[[T_i]][T_i^{-1}]^{\vee}\text{   for $2\leq i< n$}   ,
   \end{align}
   and define
   \begin{align}
    \fOK &= B_{n-1},\\
    \fK  &= \fOK [p^{-1}],
   \end{align}
   where $^{\vee}$ denotes the $p$-adic completion. So $\fK$ is a complete $n$-dimensional local field
   with residue field 
   \begin{equation}
    k= \fOK/(p) \cong \FF_p((T_1))\dots ((T_{n-1})).
   \end{equation}
   
   \noindent {\bf Remark.} The topology on $\fK$ is given by the $p$-adic valuation $v$.
   \vs
   
   We also define
   \begin{align}
    A_1 &= \ZZ_p[[T_1]][T_1^{-1}]^{\vee},\\
    A_i &= A_{i-1}[[T_i]][T_i^{-1}]^{\vee}\text{   for $2\leq i< n$}   ,\\
    \fOH &= A_{n-1},\\
    \fH  &= \fOH [p^{-1}],
   \end{align}
   so $\fK$ is a finite extension of $\fH$ of degree $[K:\Qp]$. Note that the residue field $h$ of $\fH$
   is isomorphic to $k$. Define the following towers of fields 
   \begin{align}
    H_k&=\Qp(\muk),\\
    K_k&=K(\muk),\\
    \fH_k&= \fH(\muk,T_1^{\frac{1}{p^k}},\dots,T_n^{\frac{1}{p^k}}),\\
    \fK_k&= \fK(\muk,T_1^{\frac{1}{p^k}},\dots,T_n^{\frac{1}{p^k}}),
   \end{align}
   and denote by $H_\infty$,$K_\infty$, $\fH_\infty$ and $\fK_\infty$ the respective direct limits. Note
   that $\fK_\infty$ is a finite extension of $\fH_\infty$ of degree $[K_\infty:H_\infty]$. For
   simplicity, we make the follwing assumption:-
   \vspace{1ex}
   
   \noindent {\bf Assumption.} The fields $H_\infty$ and $K_\infty$ have the same residue field 
   (which is isomorphic to $\Fp((T_1))\dots((T_{n-1}))$).
   \vspace{1ex}
   
   Let $\bar{\pi}_K=(\bar{\pi}_{K,n})_{n\geq 1}$ be a norm-compatible system of uniformizers of the tower
   $(K_n)_{n\geq 1}$. 
   
   \begin{lem}
    The field of norms $\EfK$ of $\fK_\infty$ is given by
    \begin{align}
      \EfK &\cong (\Fp[[\bar{\pi}_K]]\hat{\otimes}\Fp((T_1))\dots((T_{n-1})))[\bar{\pi}_K^{-1}]\\
           &\cong \Fp((T_1))\dots((T_{n-1}))((\bar{\pi}_K)).
    \end{align}
   \end{lem}
   \begin{proof}
    See the section on Kummer towers in~\cite{scholl}. (Note
    that $\Fp((T_1))\dots((T_{n-1}))$ has the discrete topology.)
   \end{proof}
   
   Similarly, the field of norms of $\fH$ is given by 
   \begin{equation}
    \EfH\cong \Fp((T_1))\dots((T_{n-1}))((\bar{\pi})),
   \end{equation}
   where $\bar{\pi}=\epsilon -1$ is the canonical system of uniformizers of the tower $(H_n)_{n\geq 1}$.
   (For the definition of $\epsilon$, see e.g.~\cite{cherbonniercolmez1}.) 
   
   We also define
   \begin{align}
    \EfK^+&=\Fp((T_1))\dots((T_{n-1}))[[\bar{\pi}_K]],\\
    \EfH^+&=\Fp((T_1))\dots((T_{n-1}))[[\bar{\pi}]].
   \end{align}

%+++++++++++++++++++++  
  
  \subsection{Lift to characteristic $0$}
  
   In this subsection, we define the higher dimensional analogues of Fontaine's rings $\AA_{\Qp}$ and
   $\AA_K$. We start with the lift of $\EfH$:-
   \vs
   
   \noindent {\bf Definition.}
   Let $\pi=[\epsilon]-1$ be the usual lift of $\bar{\pi}$ to characteristic $0$ (where $[\epsilon]$
   denotes the Teichmueller representative of $\epsilon$). Let $\ffR^+=\OfH[[\pi]]$ and
   $\ffR=\ffR^+[\pi^{-1}]^\vee$, so $\ffR$ is a complete regular local ring of dimension $n$ with residue
   field $\ffR/(p)\cong \EfH$. 
   \vs
   
   \noindent {\bf Note.} $\AA_{\Qp}\subset \ffR$ and $\AA_{\Qp}^+\subset\ffR^+$.
   \vs

%++++++++Action of Frobenius++++++++++   
   
   \noindent {\bf Action of Frobenius.} The action of the Frobenius operator $\phi$ on $\ffR$ is
   determined by
   \begin{align}
    \phi: T_i&\rightarrow T_i^p\text{   for all $i$}\\
          \pi&\rightarrow (\pi+1)^p-1.
   \end{align}
   
   \begin{lem}
    The ring $\ffR$ is a free $\phi(\ffR)$-module of rank $p^n$, and a basis is given by 
    $\{ T_1^{i_1}\dots T_{n-1}^{i_{n-1}}(\pi+1)^j\}_{0\leq i_1,\dots,i_{n-1},j\leq p-1}$.
   \end{lem}
   \begin{proof}
    We can assume that the action of $\phi$ on $T_i$ is trivial for all $i$. The lemma then follows by
    explicit calculation.
   \end{proof}
   
%++++++++Lift of E_K+++++++++

   \noindent {\bf Construction of the lift $\fR_K$ of $\EfK$.} We start by proving the following lemma:-
   
   \begin{lem}
    The field $\EfK$ is given by
    \begin{align}
     \EfK&\cong \EK\otimes_{\EH}\EfH\\
         &= \Fp((\bpi_K))\otimes_{\Fp((\bpi))}\Fp((X_1))\dots((X_{n-1}))((\bpi))
    \end{align}
   \end{lem}
   \begin{proof}
    We have $[\EfK:\EfH]=[\EK:\EH]=[K_\infty:H_\infty]$, and there is a natural injection
    \begin{equation}
     \Fp((\bpi_K))\otimes_{\Fp((\bpi))}\Fp((X_1))\dots((X_{n-1}))((\bpi))\hookrightarrow 
     \Fp((X_1))\dots((X_{n-1}))((\bpi_K)),
    \end{equation}
    given by `component-wise' multiplication.
%    (regarding the elements of $\Fp((X_1))\dots((X_{n-1}))((\bpi))$ as Laurent 
%    series in $T$). 
    But $\Fp((\bpi_K))$ is an $\Fp((\bpi))$-vector space of dimension 
    $[K_\infty:H_\infty]$, so 
    \begin{equation}
     \dim_{\EfH}\Fp((\bpi))\otimes_{\Fp((\bpi))}\EfH=[K_\infty:H_\infty].
    \end{equation}
   \end{proof}
   
   We can now make the following definitions:-
   \vs
   
   \noindent {\bf Definition.} Let $\fR_K=\AA_K\otimes_{\AA_{\Qp}}\ffR$, and let $\fS_K=\fR_K[p^{-1}]$.
   \vs

   Then $\fS_K$ is an $n+1$-dimensional local field (in the $p$-adic topology) and
   $\fR_K$ is a $p$-adically complete local ring of dimension $n+1$ and 
   residue field $\EfK$, so it is a lift of $\EfK$ to characteristic $0$. 
   The natural topology on $\fR_K$ is weaker than the 
   $p$-adic toplogy:- A family of neighbourhoods of $0$ in $\fR_K$ is given by
   $p^k\fR_K+ \pi_K^lO_{\fH}[[\pi_K]]$ for $k,l\geq 1$. (On can show that this definition is
   independent of the lift $\pi_K$ of $\bpi_K$.) The natural topology on $\fS_K$ is then the inductive
   limit topology.
   \vs
   
   The action of the Frobenius operator $\phi$
   extends naturally to $\fR_K$ and $\fS_K$, with the usual action on $\AA_K$. 
   Note that the action of $\phi$ restricts to $\fR_K^+=\{ a\in\fR_K: a\mod p\in\EfK^+\}$.

   \begin{lem}\label{phimodule+}
    $\fR_K^+$ is a free $\phi(\fR_K^+)$-module of rank $p^n$.
   \end{lem}
   \begin{proof}
    It follows from the classical theory that 
    $\AA_K^+$ is a free $\phi(\AA_K^+)$-module of dimension $p$ 
    with basis $1,\pi+1,\dots,(\pi+1)^{p-1}$, and
    $\ffR^+$ is a free $\phi(\ffR^+)$-module of dimension $p^n$
    with basis 
    $\{ T_1^{i_1}\dots T_{n-1}^{i_{n-1}}(\pi+1)^j\}_{0\leq i_1,\dots,i_{n-1},j\leq p-1}$.
   \end{proof}
   
   Similarly, we have
   
   \begin{lem}\label{phimodule}
    $\fR_K$ is a free $\phi(\fR_K)$-module of rank $p^n$, and $\fS_K$ is a $\phi(\fS_K)$-vector space of
    dimension $p^n$.
   \end{lem}
   
%+++++++++++++++++++++++++++
   
  \subsection{Some details about $\phi(\pi_K)$}\label{detailsFrobenius}
  
   Recall the following result from~\cite{berger1}:-
    
   \begin{prop}\label{quoteberger1}
    Let $e=[K:\Qp]$. Then
    there exists $n(K)\in\NN$ and an element $\pi_K\in\AA_K^{\dagger,n(K)}$ such
    that the image of $\pi_K$ mod $p$ is the uniformizing element $\bpi_K$ of $\EE_K$. 
    Furthermore, if $n\geq n(K)$, then $\AA_K^{\dagger,n}$ is isomorphic to the ring of series 
    $\sum_{k\in\ZZ}c_kX^k$, $c_k\in\Zp$ for all $k$, which are analytic and bounded above by $1$ on 
    the annulus $\{ p^{-\frac{1}{e(p-1)p^{n-1}}}\leq\mid X\mid<1\}$.
   \end{prop}
   
   \noindent {\bf Remark.} This result is hidden in the proof of Proposition III.2.1
   in~\cite{cherbonniercolmez2} - the main point is that for $n\gg 0$, we have 
   $\AA^{\dagger,n}_K=\AA^{\dagger,n}_{\QQ_p}[\pi_K]$, where $\pi_K\in\AA_K$ is a lift of $\bpi_K$. 
   \vs
   
   Proposition~\ref{quoteberger1} has the following useful consequence:-
   
   \begin{cor}
    We have that $\AA_K=\AA_{\Qp}[\pi_K]$. It follows that any element of
    $\AA_K$ - and hence any element of $\fR_K$ - can be written as a power series in $\pi_K$.
   \end{cor}
   
   It is shown in~\cite{cherbonniercolmez1} that the power series from Propsition~\ref{quoteberger1} 
   can be characterized as follows:-
   
   \begin{lem}
    Let $x=\sum_{k\in\ZZ}c_kX^k$ with $c_k\in\Zp$ for all $k$. Then $x$ is analytic and bounded above by
    $1$ on 
    the annulus $\{ p^{-\frac{1}{e(p-1)p^{n-1}}}\leq\mid X\mid<1\}$ if and only if
    $v_p(c_k)+\frac{k}{(p-1)ep^{n-1}}\geq 0$ for all $k<0$ and $\rightarrow +\infty$ as $k\rightarrow
    -\infty$.
   \end{lem}
   
   \begin{cor}\label{invertibility}
    Let $n\geq n(K)$, and let $x=\sum_{k\in\ZZ}c_k\pi_K^k\in\AA_K^{\dagger,n}$. If $c_0\in\ZZ_p^\times$, 
    then $x$ is invertible in $\AA_K^{\dagger,n+1}$.
   \end{cor}
   
   For the rest of this section we fix this lift $\pi_K$ and assume that $n>n(K)$. Our
   aim is to prove Corollary~\ref{structurederivative} below, which gives some details
   about the action of $\phi$ on $\pi_K$. Recall that by our
   simplifying assumption $e=[K_\infty:F_\infty]$.
   \vs 
   
   To prove Corollary~\ref{structurederivative}, we first need a couple of lemmas:-
   
   \begin{lem}\label{formalderivative}
    Let $n\gg 0$, and let $x\in\AA_K^{\dagger,n}$. Write $x$ as a power series in $\pi_K$, i.e.
    \begin{equation}
     x=\dots+c_{-2}\pi_K^{-2}+c_{-1}\pi_K^{-1}+c_0+c_1\pi_K+\dots,
    \end{equation}
    and let 
    \begin{equation}
     x'=\dots-2c_{-2}\pi_K^{-3}-c_{-1}\pi_K^{-2}+c_1+2c_2\pi_K+\dots
    \end{equation}
    be the formal derivative $x'=\frac{dx}{d\pi_K}$. Then $x'\in\AA_K^{\dagger,n}$.
   \end{lem}
   \begin{proof}
    By the definition of $\AA_K^{\dagger,n}$, we know that
    \begin{align}
     v_p(c_i)+\frac{i}{(p-1)ep^{n-1}}\geq 0 &\text{   for all $i<0$}\\ 
     v_p(c_i)+\frac{i}{(p-1)ep^{n-1}}\rightarrow +\infty &\text{   as $i\rightarrow -\infty$}
    \end{align}
    Now $v_p(c_i)$ is always a positive integer, so if $((p-1)ep^{n-1})\nmid i$, then 
    \begin{equation}
     v_p(c_i)+\frac{i}{(p-1)ep^{n-1}}\geq \frac{1}{(p-1)ep^{n-1}}>0,
    \end{equation}
    so certainly 
    \begin{equation}
     v_p(i c_i)+\frac{i-1}{(p-1)ep^{n-1}}\geq 0.
    \end{equation}
    Now if $((p-1)ep^{n-1})\mid i$, then $v_p(i)\geq n-1$, so
    \begin{equation}
     v_p(i c_i)+\frac{i-1}{(p-1)ep^{n-1}}\geq 0.
    \end{equation}
    Also, since 
    \begin{equation}
     v_p(c_i)+\frac{i}{(p-1)ep^{n-1}}\rightarrow +\infty \text{   as $i\rightarrow -\infty$},
    \end{equation}
    we see that
    \begin{equation}
     v_p(c_i)+\frac{i-1}{(p-1)ep^{n-1}}\rightarrow +\infty \text{   as $i\rightarrow -\infty$},
    \end{equation}
    so certainly
    \begin{equation}
     v_p(i c_i)+\frac{i-1}{(p-1)ep^{n-1}}\rightarrow +\infty \text{   as $i\rightarrow -\infty$},
    \end{equation}
    which finishes the proof.
   \end{proof}
   
   \begin{lem}\label{dividebyp}
    Let $x\in\AA_K^{\dagger,N}$ for some $N\gg 0$. Write 
     $$x=\dots+a_{-1}\pi_K^{-1}+a_0+a_1\pi_K+a_2\pi_K^2+\dots,$$
    and assume that $x$ is precisely divisible by $p^k$ in $\AA_K$. 
    Then there exists $M>0$ and $L\geq 0$ such that 
    $$x=p^k\pi_K^{-L}y$$ for some $y\in(\AA_K^{\dagger,N+M})^\times$.
   \end{lem}
   \begin{proof}
    Let $L$ be minimal such that $v_p(a_L)=k$. Suppose first that $L=0$. Now
    $$v_p(a_i)-k+\frac{i}{(p-1)ep^{N-1}}\rightarrow +\infty\text{   when $i\rightarrow -\infty$},$$
    so there exists $K\gg 0$ such that $v_p(a_i)-k+\frac{i}{(p-1)ep^{N-1}}\geq 0$ for all $i\leq -K$.
    Choose $M>0$ such that
    $$\frac{K}{(p-1)ep^{N+M-1}}<1.$$
    We can procede analogously if $L<0$. If $L>0$, then choose $M$ such that $p^M>L+1$.
    Since by the choice of $L$ the constant term of the $\pi_K$-expansion will be in $\ZZ_p^\times$, it
    follows from Corollary~\ref{invertibility} that $y$
    will be invertible in $\AA_K^{\dagger,N+M}$ if $M$ is sufficiently large.
   \end{proof}
    
   We assume for now that $p\mid e$. (We can argue similarly when $e$ and
   $p$ are coprime.)
   By the theory of ($1$-dimensional)
   fields of norms, $\EE_K$ is a separable extension of $\EE_{\Qp}$ of degree $e$, so since $\bpi_K$
   is a uniformizer of $\EE_K$, it satisfies a separable polynomial
   \begin{equation}\label{minimalreduction}
    \bar{f}(X)=X^e+\bar{a}_{e-1}X^{e-1}+\dots+\bar{a}_0,
   \end{equation}
   where $\bar{a}_i\in\EE_{\Qp}^+$ for all $i$. Since $\bar{f}(X)$ is separable, there exists $0\leq i\leq e$
   such that $p\nmid i$ and $\bar{a}_i\neq 0$. Now $\bar{a}_0$ and $\bpi$ have the same valuation in
   $\EE_K$, so $\bar{a}_0=\bpi\times a$, where $a\in 1+\bpi_K\EE_K^+$ (or maybe multiplied by an element in
   $\FF_p^\times$, but we ignore that). Using~\eqref{minimalreduction}, one can now expand $\bpi$ as a power
   series in $\bpi_K$,
   \begin{equation}\label{expansionmodp}
    \bpi=\bpi_K^e(1+\bar{b}_{1}\bpi_K+\bar{b}_{2}\bpi_K^{2}+\dots),
   \end{equation}
   where $\bar{b}_i\in\FF_p$. By the remark following~\eqref{minimalreduction}, there exists $i>0$ such
   that $\bar{b}_{i}\neq 0$ and $p\nmid i$. Let $I$ be the minimal such integer.
   \vs
   
   Now $\pi_K\in\AA_K$ is a lift of $\bpi_K$, so it satisfies a polynomial
   \begin{equation}\label{minimal}
    f(X)=X^e+a_{e-1}X^{e-1}+\dots+a_0,
   \end{equation}
   where $a_i\in\AA_{\Qp}$ and $f(X)\mod p=\bar{f}(X)$. In fact, it is shown in the proof of
   Proposition III.2.1 in~\cite{cherbonniercolmez2} that $a_i\in\AA_{\Qp}^+$. 
   As above, we can use equation~\eqref{minimal} to
   expand $\pi$ as a power series in $\pi_K$,
   \begin{equation}\label{expansion}
    \pi=\pi_K^e(b_0+b_1\pi_K+b_2\pi_K^2+\dots),
   \end{equation}
   which reduces to~\eqref{expansionmodp} when taken $\mod p$, i.e. 
   \begin{align}
    &b_0=1     \mod p,\\ 
    &b_i\mod p=\bar{b}_i &\text{for all $i\geq 1$}.
   \end{align}
   Applying $\phi$ to~\eqref{expansion} gives
   \begin{equation}\label{phiexpansion}
    (\pi+1)^p-1=\phi(\pi_K)^e(b_0+b_1\phi(\pi_K)+b_2\phi(\pi_K)^2+\dots).
   \end{equation}
   Write 
    $$\phi(\pi_K)=\dots+c_{-1}\pi_K^{-1}+c_0+c_1\pi_K+c_2\pi_K^2+\dots$$ 
   with $c_i\in\ZZ_p$. (Such an expansion exists since $\phi(\pi_K)\in\AA_K^{\dagger,n}$ for all 
   $n\geq 2$, and by~\ref{quoteberger1} all elements of $\AA_K^{\dagger,n}$ can be written as power 
   series in $\pi_K$ for all $n\gg 0$.)
   
   \begin{prop}\label{notpdivisible}
    There exists $i\leq -1$ such that $p\nmid i$, $p\mid c_i$ but $p^2\nmid c_i$.
   \end{prop}
   \begin{proof}
    Recall that $I$ is the
    minimal positive integer not divisible by $p$ such that $p\nmid b_I$. Expanding the left hand
    side of~\eqref{phiexpansion} as a power series in $\pi_K$ shows that the coefficient of $\pi_K^I$ is
    precisely divisible by $p$. 
    Note that $\phi(\pi_K)\mod p=\bpi_K^p$, i.e. $p\mid c_i$ for all $i\neq p$. Since $p\mid e$, it
    follows that in the power series expansion of $\phi(\pi_K)^e$, every coefficient apart from the one of
    $\pi_K^{pe}$ is divisible by $p^2$. Now if $p\mid b_j$, then it is clear that the
    coefficient of $\pi_K^I$ in $b_j\phi(\pi_K)^{e+j}$ is divisible by $p^2$. (Recall that $p\mid c_i$ for
    all $i\neq p$.) Similarly, if $j$ is a positive integer divisible by $p$, then the coefficient of 
    $\pi_K^I$ in $b_j\phi(\pi_K)^{e+j}$ is divisible by $p^2$. Recall that by the choice of $I$, $p\mid
    b_j$ for all $0<j<I$. 
    Since in the power series expansion of $\phi(\pi_K)^e$ every coefficient apart from the one of
    $\pi_K^{pe}$ is divisible by $p^2$, there must exist some $j\geq I$ not divisible by $p$ such that
    $c_{-p(e+j-1)+e+I}$ is precisely divisible by $p$. Now by assumption $e\geq p$ and $p\neq 2$, 
    which implies that $p(e+j-1)>e+I$ and finishes the proof.
   \end{proof}
   
   \begin{cor}\label{structurederivative}
    Let $\phi(\pi_K)=\dots+c_{-1}\pi_K^{-1}+c_0+c_1\pi_K+c_2\pi_K^2+\dots$, where $c_i\in\Zp$ for all $i$, 
    and let
    $g(\pi_K)=\frac{\text{d}}{\text{d}\pi_K}\phi(\pi_K)$ be its formal derivative. 
    Then there exists $I<0$ such that
    \begin{equation}
     g(\pi_K)=p\pi_K^{I}\times a,
    \end{equation}
    where $p\nmid I$ and $a\in(\AA_K^{\dagger,n})^\times$ for some $n\gg 0$.
   \end{cor}
   \begin{proof}
    Let $I<0$ be minimal such that $p\nmid I$ and $c_I$ is precisely divisible by $p$ (i.e. not divisible
    by $p^2$). Then we can write
    \begin{equation}
     g(\pi_K)=p\pi_K^Ia,
    \end{equation}
    and since $p\nmid I$, we have $a\in(\AA_K^{\dagger,n})^\times$ for $n\gg 0$ by Lemma~\ref{dividebyp}.  
   \end{proof}
    
%+++++++++++++++++++++++++++++++++++++++++++++++++++++++++++++++++++++++++++++++++++++++++++
%+++++++++++++++++++++++++++++++++++++++++++++++++++++++++++++++++++++++++++++++++++++++++++
 
  \section{Some stuff we need}
 
%++++++++++++++++++++++++
 
   \subsection{Milnor $K$-groups and completions}
 
    \noindent {\bf Notation.} For any commutative 
    ring $A$, denote by $K_n(A)$ Milnor's $K_n$-group of $A$.
    \vs
    
    Recall that $K_n(A)$ is the group generated by symbols $\{ a_1,\dots,a_n\}$, $a_i\in A^\times$ for all
    $i$, subject to the relations
    
    \noindent (1) $\{ a_1,\dots,a_n\}$ is multiplicative in each $a_i$;
    
    \noindent (2) $\{ a_1,\dots,a_n\}=1$ if $a_i+a_j=1$ for some $i\neq j$.
    
    It follows from these relations that interchanging any two entries of $a=\{ a_1,\dots,a_n\}$ yields the
    inverse $a^{-1}$. In particular, if $a_i=a_j$ for some $i\neq j$, then $a^2=1$.
    \vs
    
    \noindent {\bf Definition.} Let $U$ be the subgroup of $\fR_K$ consisting of the elements $x$ such
    that $x\mod p\in\EfK^+$, and let $\fU$ be the subgroup of $K_n(\fS_K)$ generated by the symbols
    $\{a_1,\dots,a_n\}$ with $a_1\in U^\times$, $a_i\in\fS_K^\times$. Define
    \begin{equation}
     \tK_n(\fS_K)=\varprojlim K_n(\fS_K)\slash p^m\fU.
    \end{equation}
    
    \noindent {\bf Definition.} Let $\fV$ be the subgroup of $K_n(\EfK)$ generated by the symbols
    $\{a_1,\dots,a_n\}$ with $a_1\in 1+\EfK^+$, $a_i\in\EfK^\times$. Define
    \begin{equation}
     \tK_n(\EfK)=\varprojlim K_n(\EfK)\slash p^m\fV.
    \end{equation}
        
    When $A$ is a local field with maximal ideal $\mm$, we also define the $\mm$-adic completion
    $\whK_n(A)$ of $K_n(A)$ as follows:- 
    \vs
    
    \noindent {\bf Definition.} For all $r\geq 1$ let $\fU^{(r)}$ be the subgroup of $K_n(A)$
    generated by the symbols of the form $\{ a_1,\dots,a_n\}$, where $a_1\in 1+\mm^r$ and $a_i\in
    A^\times$ for all $i$. Then define
    \begin{equation}
     \whK_n(A)=\varprojlim K_n(A)\slash \fU^{(r)}.
    \end{equation}
    
    In particular, we will be interested in the case when $A=\EfK$:- 
    
    \begin{lem}\label{topologiesmodp}
     We have a continuous map $$\tK_n(\EfK)\rTo\whK_n(\EfK).$$
    \end{lem}
    \begin{proof}
     Let $U^{(r)}=1+\bpi_K^r\EfK^+$, and let $\fU^{(r)}$ be the subgroup of $K_n(\EfK)$ generated by
     symbols of the form $\{ a_1,\dots,a_n\}$ where $a_1\in U^{(r)}$ and $a_i\neq 0$.
     It follows from the definitions of the filtrations that 
     for all $k\geq 1$ we have $p^kK_n(\EfK)\subset
     \fU^{(k)}$. 
    \end{proof}
    
    \noindent {\bf Note.} The above map has dense image.
    
%++++++++++++++++++++++++++
      
   \subsection{Modules of differentials}
  
    \noindent {\bf Definition.} Let $\Omega^{n-1}_{\fR_K\slash\ZZ}$ be the module of 
    absolute differential $(n-1)$-forms of $\fR_K$, and define $\OmegaRK$ to be its $p$-adic completion
    \begin{equation}
     \OmegaRK=\varprojlim\Omega^{n-1}_{\fR_K\slash\ZZ}\slash p^m\Omega^{n-1}_{\fR_K\slash\ZZ}.
    \end{equation}
          
    \begin{lem}\label{structuredifferentials}
     We have an isomorphism of $\fR_K$-modules
     \begin{equation}
      \OmegaRK\cong\bigoplus\fR_Kd\log(\pi_K)\wedge d\log(T_{i_1})\wedge\dots\wedge d\log(T_{i_{n-2}})
      \oplus\fR_K \bigwedge_{1\leq i<n}d\log(T_i)
     \end{equation}
    \end{lem}
    \begin{proof}
     Follows from explicit computation, using the observation 
     that $\pi_K+1$ is invertible in $\fR_K^+$. (Note that $d\log(\pi_K)=\frac{\pi_K+1}{\pi_K}
     d\log(\pi_K+1)$.)
    \end{proof}
         
    Now $\fR_K$ is finite flat and locally of complete intersection over
    $\phi(\fR_K)$, 
    so as shown in~\cite{fukaya1}, we have a trace map of differential forms
    \begin{equation}
     \Tr^{(n-1)}_\phi:\OmegaRK\rTo \OmegaRK
    \end{equation}
    and of course the usual trace map
    \begin{equation}
     \Tr_\phi: \fR_K\rTo \fR_K.
    \end{equation}
    \vs

    \begin{lem}\label{explicittrace}
     The trace map is characterized as follows:-
     \begin{align*}
      \Tr^{(n-1)}_\phi(\alpha.\bigwedge d\log(T_{i_j}) \wedge d\log (\pi_K+1))
      &=\frac{1}{p^{n-1}}\Tr_\phi (p.\alpha\omega)\bigwedge d\log(T_{i_j})\wedge d\log (\pi_K+1),\\ 
      \Tr^{(n-1)}_\phi(\beta d\log(T_1)\wedge\dots\wedge d\log (T_{n-1}))
      &=\frac{1}{p^{n-1}}\Tr_\phi(\beta)d\log(T_1)\wedge\dots\wedge d\log (T_{n-1}),
     \end{align*}
     where
     \begin{equation}
      \omega=\frac{\phi(\pi_K)+1}{\pi_K+1}(\frac{d}{d\pi_K}\phi(\pi_K+1))^{-1}.
     \end{equation}
    \end{lem}
    \begin{proof}
     First note that since $\phi$ is a ring homomorphism, we have
     $\frac{d}{d\pi_K}\phi(\pi_K+1)=\frac{d}{d\pi_K}\phi(\pi_K)$. Now by
     Corollary~\ref{structurederivative}, there exists $I<0$ and $a\in(\AA_K^{\dagger,N})^\times$ for 
     $N\gg 0$ such that $$\frac{d}{d\pi_K}\phi(\pi_K)=p\pi_K^Ia,$$ so $p\omega\in\fS_K^{\dagger,N}$ for
     $N\gg 0$. The lemma is now an immediate consequence of Remark (iii) in Section 2.1
     in~\cite{fukaya1}.
    \end{proof}

    The following corollary will be important later:-
     
    \begin{cor}\label{contractiontrace}
     $\Tr^{(n-1)}_\phi (\OmegaRK)\subset p\OmegaRK$
    \end{cor}
    \begin{proof}
     This follows from the equations in Lemma~\ref{explicittrace} and the observation that 
     $\Tr_\phi(\OmegaRK)\subset p^n\OmegaRK$.
    \end{proof}

%++++++++++++++++++++++

   \subsection{The norm map}
   
    It follows from Milnor $K$-theory that the $p$-adic valuation on $\fS_K$ induces a surjective
    map
    \begin{equation}
     \lambda:K_n(\fS_K)\rTo K_n(\EfK)
    \end{equation}
    satisfying the following condition:- If $a_i\in\fR_K^\times$ and $\bar{a}_i$ denotes the image of
    $a_i$ in $\EfK$, then 
    \begin{equation}
     \lambda: \{a_1p^{i_1},\dots,a_np^{i_n}\}\rTo \{\bar{a}_1,\dots,\bar{a}_n\}.
    \end{equation}
    Clearly $\lambda$ extends to a map on the $p$-adic completions
    \begin{equation}
     \lambda: \tK_n(\fS_K)\rTo\tK_n(\EfK).
    \end{equation}
    
    \noindent {\bf Note.} Every element in $\tK_n(\fS_K)$ can be written as the (infinite) product of
    symbols $\{a_1,\dots,a_n\}$ where each $a_i$ either lies in $\fR_K^\times$ or is equal to $p$. 
    \vs
        
    Let $\fV$ be the subgroup of $\tK_n(\fS_K)$ topologically generated by the
    symbols $\{a_1,\dots,a_n\}$, where at least one $a_i$ is equal to $p$. Note that by the definition of
    $\lambda$, we have $\fV\subset\ker(\lambda)$. Equip $\tK_n(\fS_K)\slash\fV$ with the product topology.
    \vs
    
    Now $\fS_K$ is a finite dimensional vector space over $\phi(\fS_K)$ of dimension $p^n$ (and $\EfK$ is
    a finite dimensional vector space of $\phi(\EfK)$ of dimension $p^n$), so as shown in Corollary 7.6.3
    in~\cite{weibel} we have the following result:-
    
    \begin{lem}
     We have a commutative diagram
     \begin{diagram}
      K_n(\fS_K) & \rTo^{\lambda} & K_n(\EfK) \\
      \dTo^{N}   &                &  \dTo^{N}  \\
      K_n(\fS_K) & \rTo^{\lambda} & K_n(\EfK) \\
     \end{diagram}
     where $N$ are the respective norm maps induced from the finie field extensions $\fS_K\slash
     \phi(\fS_K)$ and $\EfK\slash \phi(\EfK)$.
    \end{lem}
    
    \noindent {\bf Note.} Since the norm maps are group homomorphisms, the above diagram extends to the
    completed $K$-groups $\tK_n(\fS_K)$ and $\tK_n(\EfK)$. 
    
    \begin{lem}
     The action of $N$ factors through $\fV$.
    \end{lem}
    \begin{proof}
     Clear from the projection formula.
    \end{proof}
    
    We want to prove the following result:-
    
    \begin{prop}\label{isomodp}
     The map $\lambda$ gives a canonical isomorphism
     \begin{equation}
      (\tK_n(\fS_K)\slash\fV)^{N=1}\cong\tK_n(\EfK).
     \end{equation}
    \end{prop}
            
    \noindent {\bf Definition.} Let $\fU$ be the subgroup of $\tK(\fS_K)\slash\fV$ topologically generated
    by the symbols of the form $\{a_1,\dots,a_n\}$, where $a_1\in 1+p\fR_K$ and $a_i\in\fR_K^\times$ for
    all $i$.
    
    \begin{lem}\label{shortexactseq}
     We have a short exact sequence
     \begin{equation}
      1\rTo\fU\rTo\tK_n(\fS_K)\slash\fV\rTo \tK_n(\EfK)\rTo 1.
     \end{equation}
    \end{lem}
    \begin{proof}
     Let $x\in\tK_n(\fS_K)$. Then
     there exist $x_i\in K_n(\fS_K)$ for $i\geq 0$ such that $x$ can be represented by
     the infinite product $\prod_{i}x_i^{p^i}$. Now each of the $x_i$ can be written as a product of
     symbols $\{a^{(i)}_1,\dots,a^{(i)}_n\}$ where either one of the $a^{(i)}_j$ is equal to $p$ or all of
     the $a^{(i)}_ j$ are invertible in $\fR_K$. It follows that  
     the natural map $\tK_n(\fR_K)\rightarrow \tK_n(\fS_K)$ induces a surjection
     $\tK_n(\fR_K)\rightarrow\tK_n(\fS_K)\slash\fV$, and we have a commutative diagram
     \begin{diagram}
      1 & \rTo & \fU & \rTo & \tK_n(\fS_K)\slash\fV & \rTo & \tK_n(\EfK)  & \rTo & 1\\
        &      & \uTo^{\alpha}&      & \uTo^{\beta} &      & \uTo^{\gamma}&      & \\
      1 & \rTo & \fW & \rTo & \tK_n(\fR_K)          & \rTo & \tK_n(\EfK)  & \rTo & 1
     \end{diagram}
     Now $\gamma$ is the identity, so $\alpha$ is surjective by the snake lemma. It is therefore
     sufficient to show that $\fW$ is topologically 
     generated by symbols of the form $\{a_1,\dots,a_n\}$, where $a_1\in 1+p\fR_K$ and 
     $a_i\in\fR_K^\times$ for all $i$. Now we have a short exact sequence
     \begin{equation}\label{noncomplete}
      1\rTo W\rTo K_n(\fR_K)\rTo K_n(\EfK)\rTo 1,
     \end{equation}
     and imitating the proof of Proposition 4.1 in~\cite{fukaya2} shows that $W$ is the subgroup generated
     by symbols of the form $\{a_1,\dots,a_n\}$, where $a_1\in 1+p\fR_K$ and 
     $a_i\in\fR_K^\times$ for all $i$. Passing to the completions of~\eqref{noncomplete} gives the
     bottom row of the above diagram and hence finishes the proof.
    \end{proof}
    
    \begin{lem}\label{trivialaction}
     The action of $N$ on $\tK_n(\EfK)$ is trivial.
    \end{lem}
    \begin{proof}
     As $\phi:\EfK\rightarrow\EfK$ is the $p$th power map, the result follows by applying Lemma 12,
     paragraph 3.3 in~\cite{kato1}.
    \end{proof}
    
    In view of the preceding two lemmas, to prove Proposition~\ref{isomodp} it is sufficient to prove the
    following result:-
    
    \begin{lem}\label{nilpotentaction}
     The action of $N$ on $\fU$ is nilpotent.
    \end{lem}
    \begin{proof}
     Let 
     \begin{equation}
      \exp_{p,n-1}:\OmegaRK\rTo\fU
     \end{equation}
     be the exponential homomorphism for the $K_n$-group defined by Kurihara in~\cite{kurihara2} which is
     characterised as follows:- For $\alpha\in\fR_K$, $\beta_1,\dots,\beta_{n-1}\in\fR_K^\times$,
     \begin{equation}
      \exp_{p,n-1}(\alpha . d\log(\beta_1)\wedge \dots\wedge d\log(\beta_{n-1}))
      =\{\exp(p\alpha),\beta_1,\dots,\beta_{n-1}\}
     \end{equation}
     Since $\OmegaRK$ is $p$-adically complete, one can check that the map is surjective. (Note that
     $\exp_{p,n-1}$ is continuous.) The equations in
     Lemma~\ref{explicittrace} show that we have a commutative diagram
     \begin{diagram}
      \exp_{p,n-1}:   &   \OmegaRK              & \rTo  &\tK_n(\fS_K)\slash\fV \\
                      &   \dTo_{\Tr^{(n-1)}_\phi} &       & \dTo_{N} \\
      \exp_{p,n-1}:   &   \OmegaRK              & \rTo  &\tK_n(\fS_K)\slash\fV \\
     \end{diagram}
     Now Corollary~\ref{contractiontrace} shows that $\Tr^{(n-1)}_\phi$ is nilpotent on $\OmegaRK$, which
     finishes the proof.
    \end{proof}
     
%++++++++++++++++++++++++++++++++++++++++++++++++

  \section{The logarithmic derivative}
 
   In this section, we will use Proposition~\ref{isomodp} to construct the (generalized) logarithmic
   derivative. But before we do this, we need to introduce another module of differentials:-
 
   \subsection{Some more modules of differentials}
   
    \noindent {\bf Definition.} Let $\Omega^{1}_{\fR_K\slash\ZZ}$ be the the module of 
    absolute differential $1$-forms of $\fR_K$, and let
    \begin{equation}
     \Omega^1_{\fR_K^+\slash\ZZ}(\log)=(\Omega^1_{\fR_K^+\slash\ZZ}\oplus\fR_K^+\otimes_{\ZZ}\fR_K^\times)
     \slash\fN,
    \end{equation}
    where $\fN$ is the $\fR_K^+$-submodule of the direct sum generated by the elements $(-da,a\otimes a)$
    for $a\in\fR_K^+\cap\fR_K^\times$. In $\Omega^1_{\fR_K^+\slash\ZZ}(\log)$, denote the 
    class $(0,1\otimes a)$ for $a\in\fS_K^\times$ by $d\log(a)$. Let
    $\OmegaRKnlog$ be its $p$-adic completion
    \begin{equation}
     \OmegaRKnlog=\varprojlim\Omega^{n}_{\fR_K^+\slash\ZZ}(\log)\slash p^m\Omega^{n}_{\fR_K^+\slash\ZZ}(\log).
    \end{equation}
    Similarly, define
    \begin{equation}
     \Omega^1_{\fR_K^+\slash\ZZ}(\log)(\log)=(\Omega^1_{\fR_K^+\slash\ZZ}\oplus\fR_K^+
     \otimes_{\ZZ}\fS_K^\times)\slash\fM,
    \end{equation}
    where $\fM$ is the $\fR_K^+$-submodule of the direct sum generated by the elements $(-da,a\otimes a)$
    for $a\in\fR_K^+\cap\fS_K^\times$. In $\Omega^1_{\fR_K^+\slash\ZZ}(\log)(\log)$, denote the 
    class $(0,1\otimes a)$ for $a\in\fS_K^\times$ by $d\log(a)$. Let
    $\OmegaRKnloglog$ be its $p$-adic completion
    \begin{equation}
     \OmegaRKnloglog=\varprojlim\Omega^{n}_{\fR_K^+\slash\ZZ}(\log)(\log)
     \slash p^m\Omega^{n}_{\fR_K^+\slash\ZZ}(\log)(\log).
    \end{equation}
    Finally, let 
    \begin{equation}
     \OmegaRKn=\varprojlim\Omega^{n}_{\fR_K^+\slash\ZZ}\slash p^m\Omega^{n}_{\fR_K^+\slash\ZZ}.
    \end{equation}
              
    \begin{lem}\label{structuredifferentials2}
     We have isomorphisms of $\fR_K^+$-modules
     \begin{align}
      \OmegaRKn& \cong \fR_K^+ d\log(\pi_K+1)\wedge d\log(T_{i_1})\wedge\dots\wedge d\log(T_{i_{n-1}}), \\
      \OmegaRKnlog&\cong \fR_K^+ d\log(\pi_K)\wedge d\log(T_{i_1})\wedge\dots\wedge d\log(T_{i_{n-1}}).
     \end{align}
    \end{lem}
    \begin{proof}
     Follows from explicit computation, using the observation 
     that $\pi_K+1$ is invertible in $\fR_K^+$. (Note that $d\log(\pi_K)=\frac{\pi_K+1}{\pi_K}
     d\log(\pi_K+1)$.)
    \end{proof}
    
    \begin{lem}
     The natural map 
     $\Omega^{n}_{\fR_K^+\slash\ZZ}(\log)\rightarrow \Omega^n_{\fR_K^+\slash\ZZ}(\log)(\log)$ is injective,
     and its cokernel is killed by $p$.
    \end{lem}
    \begin{proof}
     Let $a\in\fR_K^+\cap\fS_K^\times$, and write $a=p^kb$, where $v(a)=0$. Then
     \begin{equation}
      (-da,a\otimes a)=p^n(-db,b\otimes b)+b(0,p^n\otimes p^n).
     \end{equation}
     Now $(0,p^n\otimes p^n)=-dp^n=0$ and $p^n(-db,b\otimes b)\in\fN$, so $(-da,a\otimes a)$ 
     is already trvial in
     $\OmegaRKn$. It follows that the above map is injective. 
     Since $p(x,y\otimes p)=0$ for all $x,y$, it follows that the cokernel of the map is killed by $p$.
    \end{proof}
    
    \begin{cor}\label{isodifferentials}
     Let $U$ be the subgroup of $\OmegaRKnloglog$ generated by the elements of the form 
     $\beta .d\log(p)\wedge d\log(b_2)\wedge\dots\wedge d\log (b_n)$.
     We have a canonical isomorphism of $\fR_K^+$-modules
     \begin{equation}
      \OmegaRKnloglog\slash U\cong\OmegaRKnlog.
     \end{equation}
    \end{cor} 
    
    We use this isomorphism to identify $\OmegaRKnloglog\slash U$ with
    $\OmegaRKnlog$.

%+++++++++++++++++++++++++
    
   \subsection{The $\psi$-operator}
   
    The module $\fR_K^+$ is finite flat and locally of complete intersection over
    $\phi(\fR_K^+)$, 
    so as shown in~\cite{fukaya1}, we have a trace map of differential forms
    \begin{equation}
     \Tr^{(n)}_\phi:\OmegaRKn\rTo \OmegaRKn
    \end{equation}
    and of course the usual trace map
    \begin{equation}
     \Tr_\phi: \fR_K^+\rTo \fR_K^+.
    \end{equation}
    \vs

    \begin{lem}\label{explicitntrace}
     The trace map is characterized as follows:-
     \begin{equation}
      \Tr^{(n)}_\phi(c. \bigwedge_i d\log(T_i)\wedge d\log (\pi_K+1))
      =\frac{1}{p^n}\Tr_\phi (p.c\omega)\bigwedge_i d\log(T_i)\wedge d\log (\pi_K+1),
     \end{equation}
     where 
     \begin{equation}
      \omega=\frac{\phi(\pi_K)+1}{\pi_K+1}(\frac{d}{d\pi_K}\phi(\pi_K+1))^{-1}.
     \end{equation}
    \end{lem}
    \begin{proof}
     First note that since $\phi$ is a ring homomorphism, we have
     $\frac{d}{d\pi_K}\phi(\pi_K+1)=\frac{d}{d\pi_K}\phi(\pi_K)$. Now by
     Corollary~\ref{structurederivative}, there exists $I<0$ and $a\in(\AA_K^{\dagger,N})^\times$ for 
     $N\gg 0$ such that $$\frac{d}{d\pi_K}\phi(\pi_K)=p\pi_K^Ia,$$ so $p\omega\in\fS_K^{\dagger,N}$ for
     $N\gg 0$. The lemma is now an immediate consequence of Remark (iii) in Section 2.1
     in~\cite{fukaya1}.
    \end{proof}
    
    Note that $d\log (\pi_K)=\frac{\pi_K+1}{\pi_K}d\log (\pi_K+1)$. Using the formulae in
    Lemma~\ref{explicitntrace}, we can therefore extend the trace map to 
    $\OmegaRKnlog$:-
    \begin{equation}
     \Tr^{(n)}_\phi(\frac{c}{\pi_K}. \bigwedge_i d\log(T_i)\wedge d\log (\pi_K+1))
          =\frac{1}{p^n}\frac{\Tr_\phi (p.c\tilde{\omega})}{\pi_K}\bigwedge_i d\log(T_i)\wedge
          d\log (\pi_K+1),
    \end{equation}
    where 
    \begin{equation}
     \tilde{\omega}=\frac{\phi(\pi_K)}{\pi_K}\frac{\phi(\pi_K)+1}{\pi_K+1}
     (\frac{d}{d\pi_K}\phi(\pi_K+1))^{-1}.
    \end{equation}
        
    \noindent {\bf Definition.} Let $\psi:\OmegaRKnlog\rTo\OmegaRKnlog$ be the map $\Tr^{(n)}_\phi$.

%+++++++++++++++++++++++++
    
   \subsection{The logarithmic derivative}
        
    Since $\fR_K^+$ is $p$-adically complete, the $d\log$-operator induces a continuous map
    \begin{equation}
     d\log:\tK_n(\fS_K)\rTo \OmegaRKnlog
    \end{equation}
    which clearly factors through $\fV$, giving
    \begin{equation}
     d\log:\tK_n(\fS_K)\slash\fV\rTo \OmegaRKnlog.
    \end{equation}
    
    \begin{lem}
     We have a commutative diagram
     \begin{diagram}
      \tK_n(\fS_K)\slash\fV     &     \rTo^N      & \tK_n(\fS_K)\slash\fV \\
      \dTo_{d\log}              &                 &    \dTo_{d\log}     \\
      \OmegaRKnlog              &  \rTo^{\psi}    &  \OmegaRKnlog       \\
     \end{diagram}
    \end{lem}
    \begin{proof}
     Explicit computation.
    \end{proof}
    
    \begin{cor}\label{mapintodiff}
     We have a canonical homomorphism
     \begin{equation}
      \tau: \tK_n(\EfK)\rTo (\OmegaRKnlog)^{\psi=1}.
     \end{equation}
    \end{cor}
    \begin{proof}
     Define $\tau$ to be the composition of $d\log$ with the inverse of the isomorphism
     $(\tK_n(\fS_K)\slash\fV)^{N=1}\cong \tK_n(\EfK)$.
    \end{proof}
    
    To define the logarithmic derivative, we need to show the following result:-
    
    \begin{lem}
     The map $\tau$ extends to a map
     \begin{equation}
      \tau: \whK_n(\EfK)\rTo (\OmegaRKnlog)^{\psi=1}.
     \end{equation}
    \end{lem}
    \begin{proof}
     Let $a\in\tK_n(\EfK)$. It follows from the construction of $\tau$ that $\tau(a)$ can be determined by
     the following procedure:- Let $\alpha$ be some lift of $a$ to $\tK_n(\fS_K)$. Then
     $\psi^k(d\log(\alpha))\rightarrow \tau(a)$ as $k\rightarrow +\infty$. (One can show that 
     $\psi^{k+1}(d\log(\alpha))-\psi^k(d\log(\alpha))\in p^k\OmegaRKnlog$ for all $k$.)
     
     For $r\geq 1$, let $\fU^{(r)}$ be the subgroup of $\tK_n(\EfK)$ topologically generated by the
     symbols $\{a_1,\dots,a_n\}$ such that $a_1\in 1+\bpi_K^r\EfK^+$. Note that 
     \begin{equation}
      \whK_n(\EfK)=\varprojlim\tK_n(\EfK)\slash \fU^{(r)}.
     \end{equation}
     Suppose that $a\in\ker(\tK_n(\EfK)\rightarrow\whK_n(\EfK))$, i.e. $a\in\bigcap_r\fU^{(r)}$. Using the
     above algorithm for determing $\tau(a)$ then shows that $\tau(a)\in\bigcap_mp^m\OmegaRKnlog$ and
     hence $\tau(a)=0$ since $\OmegaRKnlog\cong\fR_K^+$. 
     Now let $x\in\whK_n(\EfK)$, and write $x$ as an (infinite) product $\prod_r x_r$, where
     $x_r\in\fU^{(r)}$. For each $r$, choose a lift $\fX_r$ in the subgroup of $\tK_n(\fS_K)$
     topologically generated by symbols of the form $\{b_1,\dots,b_n\}$ with 
     $b_1\in 1+\pi_K^rO_H[[\pi_K]]$. Let $y_r=d\log(\fX_r)$. Then
     $y_r\in\pi_K^{r}O_H[[[\pi_K]].\bigwedge d\log(T_i)\wedge d\log(\pi_K)$, so 
     the sum $\sum_r y_r$ converges to an 
     element $y\in\OmegaRKnlog$. Now the map $\Tr_\phi$ is $(p,\pi_K)$-continuous, so for all $k$ we have
     \begin{equation}
      \psi^k(y)=\psi^k(y_1)+\psi^k(y_2)+\dots
     \end{equation}
     But we have $\psi^{k+1}(y_i)-\psi^k(y_i)\in p^k\OmegaRKnlog$ for all $i$. Now $\OmegaRKnlog$ is
     $p$-adically complete, so $\psi^k(y)$ converges towards some element $z\in\OmegaRKnlog$ as
     $k\rightarrow +\infty$. Define $\tau(x)=z$.
    \end{proof}
        
    \begin{prop}\label{isoLCFT}
     We have a canonical homomorphism
     \begin{equation}
      \varprojlim\whK_n(\fK_m)\cong\whK_n(\EfK).
     \end{equation}
    \end{prop}
    \begin{proof}
     For $i\geq N$, define the evaluation map
     \begin{equation}
      h^{(N)}_i:\fTKdaggern \rTo \fK_i
     \end{equation}
     by $h^{(N)}_i(\pi_K)= \pi_{K_n}$ and $h^{(N)}_i(T_j)=T_j^{\frac{1}{p^i}}$.
     Then for all $i\gg 0$, the map $h^{(N)}_i$ induces a homomorphism
     \begin{equation}
      H^{(N)}_i:\whK_n(\fTKdaggern)\rTo \whK_n(\fK_i).
     \end{equation}
     One can show that there exists $j_i^{(n)}\geq i-n-1$ such that 
     the image of $\ker(H^{(N)}_i)$ in $\whK_n(\EfK)$ is the
     subgroup $\fA^{(N)}_i$ generated by the symbols of the form $\{ a_1,\dots,a_n\}$ with 
     $a_1\in 1+\bpi_K^{j_i^{(n)}}\EfK^+$, and we have commutative diagram
     \begin{diagram}
      K_n(\EfK)\slash\fA^{(N)}_{i+1}&\lTo &\whK_n(\fTKdaggern)\slash\ker(H^{(N)}_{i+1}) 
      &\rTo^\cong_{H^{(N)}_{i+1}}&   \whK_n(\fK_{i+1})    \\
        \dTo                        &     &                       \dTo_N                      &            &    \dTo_{N_{O_{\fK_{m+1}}\slash O_{\fK_m}}}    \\
      K_n(\EfK)\slash\fA^{(N)}_i    &\lTo &\whK_n(\fTKdaggern)\slash\ker(H^{(N)}_i)     
      &\rTo^{\cong}_{H^{(N)}_i}  &   \whK_n(\fK_i)\\
     \end{diagram} 
     Here, the map $K_n(\EfK)\slash\fA^{(N)}_{i+1}\rightarrow K_n(\EfK)\slash\fA^{(N)}_i$ is the natural
     projection map.
     Now $j_i^{(n)}\rightarrow +\infty$ as $i\rightarrow+\infty$, so we have
     $\varprojlim K_n(\EfK)\slash\fA^{(N)}_i =\whK_n(\EfK)$, and hence
     the above diagram induces a map
     \begin{equation}
      g:\varprojlim \whK_n(\fK_i)\rTo\whK_n(\EfK).
     \end{equation}
    \end{proof} 
     
    We can now define the logarithmic derivative:-
    \vs
     
    \noindent {\bf Definition.} The logarithmic derivative is the map 
    \begin{equation}
     \Tau:\varprojlim\whK_n(O_{\fK_m})\rTo (\OmegaRKnlog)^{\psi=1}
    \end{equation}
    obtained by composing the homomorphism $g$ with $\tau$.
    \vs
    
    \noindent {\bf Problem.} Recall that
    \begin{equation}
     \OmegaRKnlog\cong\fR_K^+d.\log (\pi_K)\wedge d\log(T_1)\wedge\dots\wedge d\log(T_{n-1}).
    \end{equation}
    In analogy to the construction of Cherbonnier and Colmez in the $1$-dimensional case, we would like to
    associate to an element in $\varprojlim\whK_n(O_{\fK_m})$ something like a function that is convergent
    on some annulus of the $p$-adic unit disc. However, the elements of $\fR_K^+$ (as
    series in $\pi_K$) do not necessarily converge anywhere. Consequently, we need to define a
    suitable subring, an $n$-dimensional ring of overconvergent series. We will do that in the next
    section.
    
%++++++++++++++++++++++++++++++     
     
 \section{Higher dimensional rings of overconvergent series}
 
 %+++++++++++++++++++++
 
  \subsection{Definition and properties}
  
   \subsubsection{Definition of $\fRdaggern$}
    
    As in the classical case, we define the elements of $\fRdaggern$ as certain functions defined on some
    annulus on the $p$-adic open unit disc:- Let $x\in\ffR$, and write
    $x=\sum_{i\in\ZZ}f_{i}(T_1,\dots,T_{n-1})\pi^i$. 
    \vs

    \noindent {\bf Definition.} Let $n\geq 1$. Then
    $x\in\fRdaggern$ if the series 
    $\sum_{i\in\ZZ}f_{i}(T_1,\dots,T_{n-1})\pi^i$ converges on the
    annulus $\{ p^{-\frac{1}{(p-1)p^{N-1}}}\leq \mid X\mid <1\}\subset\CC_p$ (with values in
    $\CC_p\hat{\otimes}_{\Zp}O_{\fH}$) and is bounded above by $1$.
    \vs 
    
    Alternatively, the elements of $\fRdaggern$ can be defined using growth conditions on the negative
    coefficients:-
    
    \begin{prop}\label{negativecoefficients}
     Let $x\in\ffR$, and write 
     $x=\sum_{i\in\ZZ}f_{i}(T_1,\dots,T_{n-1})\pi^i$. 
     Then $x\in\fRdaggern$ if and only if 
     $v(f_i)+\frac{i}{(p-1)p^{N-1}}\geq 0$ for all $i<0$ and $\rightarrow +\infty$ as 
     $i\rightarrow -\infty$.
    \end{prop}   
    \begin{proof}
     Explicit calculation.
    \end{proof}   
    
    \begin{cor}\label{moduleFontaine}
     The ring $\fRdaggern$ is an $\AA^{\dagger,N}_{\Qp}$-module.
    \end{cor}

    It follows from the definition and the properties of the ring $\AA_{\Qp}^{\dagger,N}$ that we have 
    the following corollary:-
   
    \begin{cor}\label{properties1}
     (1) $\fRdaggern\subset \ffR^{\dagger, N+1}$, and
   
     \noindent (2) If $x\in\fRdaggern$, then $\phi(x)\in\ffR^{\dagger,N+1}$.
    \end{cor}
    \begin{proof}
     To prove (2), observe that if $\alpha\in \{ p^{-\frac{1}{(p-1)p^N}}\leq \mid X\mid <1\}$, then 
     $(\alpha+1)^p-1\in \{ p^{-\frac{1}{(p-1)p^{N-1}}}\leq \mid X\mid <1\}$.
    \end{proof}
    
    \begin{prop}\label{moduleoverphi1}
     The ring $\ffR^{\dagger,N+1}$ is a free $\phi(\fRdaggern)$-module of degree $p^n$.
    \end{prop}
    \begin{proof}
     Explicit calculation shows that
     $\{ T_1^{i_1}\dots T_{n-1}^{i_{n-1}}(\pi+1)^j\}_{0\leq i_1,\dots,i_{n-1},j\leq p-1}$ 
     is a basis of $\ffR^{\dagger,N+1}$ over $\phi(\fRdaggern)$.
    \end{proof}

    \noindent {\bf Definition.} The ring $\fRdagger\subset\ffR$ of overconvergent elements is 
    the union $$\fRdagger=\bigcup\fRdaggern.$$

    The topology on $\fRdaggern$ is induced from the topology on $\ffR$:-
    
    \begin{lem}
     A set of neighbourhoods of $0$ in $\fRdaggern$ is given by $p^k\fRdaggern+\pi^l\ffR^+$.
    \end{lem}

    \noindent {\bf Remark.} The topology on $\fRdagger$ is the Frechet topology induced from the
    topologies on the $\fRdaggern$.
    \vs

%++++++++++++++

   \subsubsection{Definition of $\fRKdaggern$}
   
    Since $\AA^{\dagger,N}_K$ is a free finitely generated module over $\AA^{\dagger,N}_{\Qp}$ for $N\gg
    0$, it seems natural to make the following definition:-
    \vs

    \noindent {\bf Definition.} $\fRKdaggern= \AA^{\dagger,N}_K\otimes_{\AA_{\Qp}^{\dagger,N}}\fRdaggern$
    
    \begin{lem}\label{properties2}
     (1) $\fRKdaggern\subset \fR_K^{\dagger, N+1}$,
   
     \noindent (2) if $x\in\fRKdaggern$, then $\phi(x)\in\fR_K^{\dagger,N+1}$, and
     
     \noindent (3) $\fR_K^{\dagger,N+1}$ is a free finitely generated module of degree $p^n$ over 
     $\phi(\fRKdaggern)$.
    \end{lem}
    \begin{proof}
     (1) and (2) are immediate consequences of Corollary~\ref{properties1} and the properties of the 
     ring $\AA_K^{\dagger,N}$. Now as shown in~\cite{cherbonniercolmez2}, Lemma V.1.4, the set $\{
     1,\dots,(\pi+1)^{p-1}\}$ is a basis of $\AA_K^{\dagger,N+1}$ over $\phi(\AA_K^{\dagger,N})$. Together
     with Proposition~\ref{moduleoverphi1}, this proves (3).
    \end{proof}
    
    Let $\fR_K^+=\AA_K^+\otimes_{\AA_{\Qp}^+} \ffR^+$.
    
    \begin{lem}
     A set of neighbourhoods of $0$ in $\fRKdaggern$ is given by $p^k\fRKdaggern+\pi^l\fR_K^+$ for 
     $l,k\geq 1$.
    \end{lem}
    
    If $n\gg 0$, then the elements of $\fRKdaggern$ also have a description as functions defined on some
    annulus on the $p$-adic open unit disc:-    
    Now let $x\in\fRKdaggern$, and let $n(K)$ be as defined in Proposition~\ref{quoteberger1}. We can then
    write $x$ as $\sum_{i\in\ZZ}f_i(T_1,\dots,T_{n-1})\pi_K^i$.
    \vs
   
    \begin{prop}\label{functions2}
     Let $N\geq n(K)$. Then $x\in\fRKdaggern$ if and only if the series
     $\sum_{i\in\ZZ}f_i(T_1,\dots,T_{n-1})X^i$ converges on the annulus 
     $\{ p^{-\frac{1}{e(p-1)p^{N-1}}}\leq\mid X\mid<1\}$ and is bounded above by $1$.
    \end{prop}

    To prove the proposition, we need the following result:-

    \begin{lem}
     We have that $x\in\fRKdaggern$ if and only if 
     $v(f_i)+\frac{i}{e(p-1)p^{N-1}}\geq 0$ for all $i<0$ and $\rightarrow +\infty$ as 
     $i\rightarrow -\infty$.
    \end{lem}
    \begin{proof}
     Imitate the proof of Proposition~\ref{negativecoefficients}.
    \end{proof}

    \begin{proof}
     Note that $1,\pi_K,\dots,\pi_K^{e-1}$ is a basis for $\AA_K^{\dagger,N}$ over
     $\AA_{\Qp}^{\dagger,N}$. Recall from~\eqref{expansion} that there exist $c_i\in\ZZ_p$ such that
     \begin{equation}
      \pi_K^e=\pi(c_0+c_1\pi_K+c_2\pi_K^2+\dots).
     \end{equation}
     It follows that if we write
     \begin{equation}
      \pi_K^e=a_0+a_1\pi_K+\dots+a_{e-1}\pi_K^{e-1}
     \end{equation}
     with $a_i\in\AA_{\Qp}^{\dagger,N}$, then $a_i\mod p\in\bpi_K\EE_{\Qp}^+$, i.e. as a
     function on $\{p^{-\frac{1}{(p-1)p^{N-1}}}\leq\mid X\mid<1\}$, $a_i$ is strictly less that $1$ for
     all $i$. Now the topology on $\AA_{\Qp}^{\dagger,N}$ agrees with the topology of uniform convergence
     on every compact subannulus of $\{p^{-\frac{1}{(p-1)p^{N-1}}}\leq\mid X\mid<1\}$. It follows that if
     we write 
     \begin{equation}
      \pi_K^j=a_{j,0}+a_{j,1}\pi_K+\dots+a_{j,e-1}\pi_K^{e-1}
     \end{equation}
     with $a_{j,k}\in\AA_{\Qp}^{\dagger,N}$, then $a_{j,k}\rightarrow 0$ as $j\rightarrow +\infty$ for all
     $0\leq k<e$.
    
     Suppose now that $\sum_{j\in\ZZ}f_j(T_1,\dots,T_{n-1})X^j$ converges on  
     $\{ p^{-\frac{1}{e(p-1)p^{N-1}}}\leq\mid X\mid<1\}$ and is bounded above by $1$. Write
     \begin{equation} 
      \sum_{j\in\ZZ}f_j(T_1,\dots,T_{n-1})X^j=
      \sum_{i=0}^{e-1}X^i\sum_{j\in\ZZ}f_j(T-1,\dots,T_{n-1})a_{j,i}.
     \end{equation}
     Now $v(f_j)>0$ for all $j<0$ and
     $\rightarrow +\infty$ as $j\rightarrow -\infty$ by the previous lemma, so 
     $f_j(T_1,\dots,T_{n-1})a_{i,j}$  is strictly less than $1$ on 
     $\{p^{-\frac{1}{(p-1)p^{N-1}}}\leq\mid X\mid<1\}$ and $\rightarrow 0$ as $j\rightarrow -\infty$.
     Since we have shown above that $a_{j,k}\rightarrow 0$ as $j\rightarrow +\infty$, it follows that
     the series $\sum_{j\in\ZZ}f_j(T-1,\dots,T_{n-1})a_{j,i}$ converges on 
     $\{p^{-\frac{1}{(p-1)p^{N-1}}}\leq\mid X\mid<1\}$ and is bounded above by $1$ and hence is an element
     of $\fR_{\Qp}^{\dagger,N}$.
     
     Conversely, if $x\in\fRKdaggern$, then it follows from the definition of $\fRKdaggern$ and
     Proposition~\ref{quoteberger1} that the expansion of $x$ as a series in $\pi_K$ converges on 
     the annulus $\{ p^{-\frac{1}{e(p-1)p^{N-1}}}\leq\mid X\mid<1\}$ and is bounded above by $1$.
    \end{proof}

    \noindent {\bf Definition.} Define the ring of overconvergent elements of $\fR_K$ as
    $$\fRKdagger=\bigcup_{n\geq n(K)}\fRKdaggern.$$
    \vs
  
%+++++++++++++++   
   
   \subsubsection{Properties of $\fRKdagger$}
   
    Throughout this section, we assume that $N\geq n(K)$. 
    \vs
    
    \noindent {\bf Definition.} In the $(p,\pi_K)$-adic topology on $\fRKdaggern$, a basis of
    neighbourhoods of $0$ is given by $p^k\fRKdaggern + \pi_K^lO_{\fK}[[\pi_K]]$, where $k,l\geq 1$.
    
    \begin{lem}
     The $(p,\pi)$-adic topology and the $(p,\pi_K)$-adic topology agree on $\fRKdaggern$. 
    \end{lem}
    \begin{proof}
     Recall that $\fRKdaggern=\AA^{\dagger,N}_K\otimes_{\AA_{\Qp}^{\dagger,N}}\fRdaggern$. The lemma is
     now a consequence of the fact that the $(p,\pi)$-adic topology and the $(p,\pi_K)$-adic topology agree
     on $\AA^{\dagger,N}_K$ (c.f. Section 2.6 in~\cite{berger1}.)
    \end{proof}
    
    \noindent {\bf Remark.} The $(p,\pi_K)$-adic topology is the topology of covergence on compact 
    subannuli of $\{ p^{-\frac{1}{e(p-1)p^{N-1}}}\leq\mid X\mid<1\}$.

    \begin{lem}
     The ring $\fRKdaggern$ is complete with respect to the $(p,\pi_K)$-adic topology. 
    \end{lem}
    \begin{proof}
     Clear from Proposition~\ref{functions2} (or the remark above).
    \end{proof}
    
    \begin{lem}\label{invertible2}
     Let $x\in\fRKdaggern$ and write $x=\sum f_i(T_1,\dots,T_{n-1})\pi_K^i$. 
     If $v(f_0)=0$ and 
     $v(f_i)+\frac{i}{(p-1)p^{N-1}}>0$ for all $i<0$, then $x$ is invertible in $\fRKdaggern$.
    \end{lem}
    \begin{proof}
     Suppose that $v(f_0)=0$ and $v(f_i)+\frac{i}{e(p-1)p^{N-1}}>0$ for all $i>0$. Since $O_{\fH}$ is a
     $p$-adically complete local ring, $f_0$ is invertible in $O_{\fH}$ (so wlog $=1$) and hence 
     $\sum_{i\geq 0}f_i(T)\pi_K^i$ in invertible in $\fRKdaggern$ (since the formal inverse converges), 
     so without loss of generality assume that $x=1+\sum_{j\leq -1}f_{j}(T_1,\dots,T_{n-1})
     \pi_K^{j}$. The formal 
     inverse of $x$ is given by $\sum_{k\geq 0}(1-x)^k=\sum_{i\in\ZZ}h_i(T_1,\dots,T_{n-1})\pi_K^i$. 
     We need to
     show that this sum converges in $\fRKdaggern$. Since $v(f_j)+\frac{j}{e(p-1)p^{n-1}}>0$ for all $j$,
     it is clear from the formal expansion that $v(h_j)+\frac{j}{e(p-1)p^{n-1}}>0$ for all $j<0$.
     It remains to show that $v(h_j)+\frac{j}{e(p-1)p^{n-1}}\rightarrow +\infty$ as $j\rightarrow
     -\infty$. Let $\alpha=\frac{1}{e(p-1)p^{n-1}}$. Since $v(f_i)$ is a positive integer for all $i$, 
     $v(f_i)+\frac{i}{e(p-1)p^{N-1}}\geq \alpha$ for all $i>0$. Fix $M>0$, and let $K\gg 0$ such that 
     $v(f_i)+\frac{i}{e(p-1)p^{N-1}}\geq M\alpha$ for all $i\geq K$. Expanding $\sum_{k\geq 0}(1-x)^k$
     then shows that $v(h_j)+\frac{j}{e(p-1)p^{N-1}}\geq M\alpha$ for all $j\geq MN$, which finishes the
     proof.
    \end{proof}
   
    \begin{lem}\label{invertible}
     Let $N\geq n(K)$. Let $x=\sum_{i\geq 0}f_i(T_1,\dots,T_{n-1})\pi_K^i\in\fRKdaggern$. Then the following are
     equivalent:-
     
     \noindent (i) $x$ is invertible in $\fRKdagger$,
     
     \noindent (ii) $x$ is invertible in $\fR_K^{\dagger,N+1}$,
     
     \noindent (iii) $v(f_0)=0$.
    \end{lem}
    \begin{proof}
     It follows from Lemma~\ref{invertible2} that if $v(f_0)=0$, then $x$ is invertible in 
     $\fR_K^{\dagger,N+1}$. Conversely, write $x=\sum_{i\in\ZZ}f_i(T_1,\dots,T_{n-1})\pi_K^i$, 
     and suppose that $x$ is invertible in $\fRKdagger$. 
     Since $p^{-1}\notin\fR_K$, we must have $v(f_0)=0$, i.e. $f_0$ is a unit in $O_H$ and hence
     invertible in $\fRKdaggern$. 
%     We can therefore without loss of generality assume that $f_0=1$. Then
%     the series $y=\sum_{j\geq 0} (-1)^j(x-1)^j$
%     converges in $\fR_K$ to $x^{-1}$, and if one writes $y=\sum_{j\in\ZZ}h_j(T)\pi_K^j$, then $h_j$ 
%     satisfies $v(h_j)+\frac{j}{e_K(p-1)p^{n-1}}\geq 0$ for all $j<0$. Note that $y$ 
%     is not necessarily an element of $\fRKdaggern$ since we do not know whether 
%     $v(h_j)+\frac{j}{e_K(p-1)p^{N-1}}\rightarrow +\infty$ as $j\rightarrow -\infty$. However, since
%     $\frac{1}{e_K(p-1)p^{N-1}}-\frac{1}{e_K(p-1)p^N}>0$, it follows that $y\in\fR_K^{\dagger,N+1}$.
    \end{proof}
    
    Lemma~\ref{invertible} has the following important consequence:-
    
    \begin{prop}\label{localcomplete}
     $\fRKdagger$ can be written as the (increasing) union of $p$-adically complete local rings.
    \end{prop}
    \begin{proof}
     Let $N\geq n(K)$. To simplify the notation, write (only in this proof) $R_N$ for $\fRKdaggern$.
     Recall that every element of $R_N$ can be written as a power series
     $\sum_{i\in\ZZ}f_i(T_1,\dots,T_{n-1})\pi_K^i$. Let 
     $$\pp_N=\{\sum_{i\in\ZZ}f_i(T_1,\dots,T_{n-1})\pi_K^i\in R_n \mid v(f_0)>0\}.$$
     Then $\pp_N$ is an ideal of $R_N$, and explicit calculation shows that it is prime. Note that
     $p\in\pp_N$. Let $R_{N,\pp_N}$ denote the localisation of $R_N$ at $\pp_N$, which by Lemma~\ref{invertible}
     is contained in $R_{N+1}$, and let $S_N$ be its completion with respect to $\pp_N$. Then $\pp_N$ is 
     still
     a local ring, and it it $p$-adically complete since $p\in\pp_N$. The topology defined by $\pp_N$ is
     the $(p,\pi_K)$-adic topology. To finish the proof of the proposition, note that 
     $S_N\subset R_{N+1}$ since $R_{N+1}$ is complete in the $(p,\pi_K)$-adic topology (and hence the
     topology on $S_N$ is the subspace topology induced from $R_{N+1}$), i.e.
     $\fRKdagger=\bigcup_{N\geq n(K)} S_N$.   
    \end{proof}
       
    \noindent {\bf Definition.} Write $\fSKdaggern$ for the ring $S_N$. 
    \vs
    
    \noindent {\bf Definition.} Let $\fTKdaggern=\fSKdaggern [\pi_K^{-1}].$
    \vs
    
    \noindent {\bf Note.} (1) Note that $\fTKdaggern=\fSKdaggern [\pi^{-1}].$
    
    \noindent (2) Written as power series in $\pi_K$, the elements of $\fTKdaggern$ can be
    interpreted as functions converging on $\{ p^{-\frac{1}{(p-1)e_Kp^{N-1}}}\leq \mid X\mid<1\}$ (but not
    necessarily bounded above by $1$).
    \vs
    
    \noindent {\bf Remark.} Since Frobenius of the residue field is the $p$-power map, we have
    $\phi(\fSKdaggern)\subset \fS_K^{\dagger,N+1}$ and $\phi(\fTKdaggern)\subset \fT_K^{\dagger,N+1}$.
      
 %+++++++++++++++++++++++
  
  \subsection{Properties of $\phi$}
  
   \begin{prop}
    Let $N\gg 0$. Then $\fR_K^{\dagger,N+1}$ is a free $\phi(\fRKdaggern)$-module of rank $p^N$.
   \end{prop}
   \begin{proof}
    It follows from the classical theory of the rings 
    $\AA_{\Qp}^{\dagger,N}$ and $\AA_K^{\dagger,N}$ that $\{T_1^{i_1}\dots T_{n-1}^{i_{n-1}}
    (1+\pi)^j\}_{0\leq i_1,\dots,i_{n-1},j\leq p-1}$ is a 
    basis of $\fR_K^{\dagger,N+1}$ over $\phi(\fRKdaggern)$. 
   \end{proof}
   
   \begin{cor}\label{freemoduledagger}
    $\fRKdagger$ is a free $\phi(\fRKdagger)$-module of rank $p^n$.
   \end{cor}
   
   The following proposition is the reason why the ring of overconvergent series can be used to construct
   the logarithmic derivative:-
   
   \begin{prop}\label{freeTnmodule}
    For $n\gg 0$, $\fT_K^{\dagger,N+1}$ is a free $\phi(\fTKdaggern)$-module of rank $p^n$.
   \end{prop}
   
   Noting that $\fTKdaggern=\fSKdaggern[\pi^{-1}]$, this 
   proposition is a consequence of the following result:-

   \begin{prop}\label{freeSnmodule}
    For $N\gg 0$, $\fS_K^{\dagger,N+1}$ is a free $\phi(\fSKdaggern)$-module of rank $p^n$.
   \end{prop}

   To prove the Proposition~\ref{freeSnmodule}, we begin with the following observations:-
    
   \begin{lem}\label{phipipower}
    Let $N\in\NN$. If we write $\pi^N=\sum_{0\leq i\leq p-1}\phi(a^{(N)}_i)(1+\pi)^i$, with
    $a^{(N)}_i\in\AA_{\Qp}^+$, then for all $0\leq i\leq p-1$ we have $a^{(N)}_i\rightarrow 0$ as 
    $N\rightarrow +\infty$.
   \end{lem}
   \begin{proof}
    $a^{(N)}_i=\psi(\pi^N(1+\pi)^{-i})$, and $\psi$ is continuous on $\AA_{\Qp}$.
   \end{proof}
   
   \noindent {\bf Remark.} It follows that if we expand $\pi^N$ as 
   $\sum_{0\leq i\leq p-1}\phi(b^{(N)}_i)\pi^i$, then $b^{(N)}_i\rightarrow 0$ as $N\rightarrow +\infty$.
   \vs
   
   \noindent We also need the following lemma:-
   
   \begin{lem}\label{pipowertozero}
    $\pi^i\rightarrow 0$ in $\fSKdaggern$ as $i\rightarrow +\infty$ for all $N\gg 0$.
   \end{lem}
   \begin{proof}
    Recall that $\pi=\pi_K^e(b_0+b_1\pi_K+b_2\pi_K^2+\dots)$ for some $b_i\in\ZZ_p$.
   \end{proof}
   
   \noindent We can now prove Proposition~\ref{freeSnmodule}:-    
   
   \begin{proof}
    Wlog assume that $\phi$ acts trivially on $T_1,\dots,T_{n-1}$, so $\{
    1,\pi,\dots,\pi^{p-1}\}$ is a basis for $\fR_K^{\dagger,N+1}$ over $\fRKdaggern$. 
    We first show that every element of $(\fR_K^{\dagger,n+1})_{\pp_{n+1}}$ can be written as
    $\sum_{0\leq i\leq p-1}\phi(x_i) \pi^i$ with $x_i\in\fSKdaggern$. 
    
    Let $x\in \fR_K^{\dagger,N+1}$ and write $x=\sum_{i\in\ZZ}f_i(T_1,\dots,T_{n-1})\pi_K^i$. Assume that $v(f_0)=0$. Let
    $x_0,\dots,x_{p-1}\in\fRKdaggern$ such that $x=\sum_{0\leq j<p}\phi(x_j)\pi^j$. 
    Since $v(f_0)=0$, $x_0^{-1}\in(\fRKdaggern)_{\pp_n}$, so 
    $\phi(x_0^{-1})\in(\fR_K^{\dagger,N+1})_{\pp_{N+1}}$. We can therefore
    write $x=\phi(x_0)(1+\pi\phi(x_0^{-1})z)$, where $z=\sum_{0<i\leq p-1}\phi(x_i)\pi^{i-1}$. 
    Let
    \begin{equation}\label{formalinverse}
     y=1-\pi\phi(x_0^{-1})z+\pi^2\phi(x_0^{-2})z^2-\pi^3\phi(x_0^{-3})z^3+\dots
    \end{equation}
    be the formal inverse of $1+\pi\phi(x_0^{-1})z$. Since
    $1+\pi\phi(x_0^{-1})z\in\fR_K^{\dagger,N+1}\backslash\pp_{N+1}$, we have 
    $y\in\fR_K^{\dagger,N+2}$, so there exist
    $y_0,\dots,y_{p-1}\in\fR_K^{\dagger,N+1}$ such that $y=\sum_{0\leq k\leq p-1}\phi(y_k)\pi^k$. We will
    show that in fact $y_k\in\fSKdaggern$ for all $k$. Let $0\leq k\leq p-1$.
    By expanding~\eqref{formalinverse} (and remembering that $\phi(\pi)=(\pi+1)^p-1$), we see that 
    $y_k$ is of the form
    \begin{equation}
     y_k=\sum_{l,m=0}^\infty p^l\pi^m\alpha^{(k)}_{l,m}
    \end{equation}
    for some $\alpha^{(k)}_{l,m}\in(\fRKdaggern)_{\pp_N}$. But 
    $p^l\pi^m\rightarrow 0$ in $\fSKdaggern$ as $l,m\rightarrow +\infty$ by Lemma~\ref{pipowertozero}, 
    so since $\fSKdaggern$ is complete, it follows that $y_k\in\fSKdaggern$. Since $\phi$ is a ring
    homomorhism, we have shown that every element $z\in(\fR_K^{\dagger,N+1})_{\pp_{N+1}}$ can be written
    as $\sum_{0\leq i\leq p-1}\phi(z_i)\pi^i$ with $z_i\in\fSKdaggern$. 
    
    It remains to show that this
    property extends to $\fS_K^{\dagger,N+1}$. Identify $\pp_{N+1}$ with its image in
    $(\fR_K^{\dagger,N+1})_{\pp_{N+1}}$, and assume that $x\in\pp_{N+1}^M$ for some $M\gg 0$. Write 
    $\sum_{0\leq i\leq p-1}\phi(x_i)\pi^i$ with $x_i\in\fSKdaggern$. By completeness, it is sufficient to
    show that $x_i\in\pp_n^L\fSKdaggern$ for some $L\gg 0$. Now since $x\in\pp_{N+1}^M$, there exist
    $a_1\in \fR_K^{\dagger,N+2}$, $a_2\in\fR_K^+$ and $N_1,N_2\gg 0$ such that 
    $x=p^{N_1}a_1+\pi^{N_2}a_2$. 
    Since $\phi$ is $\Zp$-linear, we can wlog. assume that $a_1=0$, so $x=\pi^Na$. But if we write 
    $\pi^N=\sum_{0\leq i\leq p-1}\phi(b^{(N)}_i)\pi^i$, then we deduce from Lemma~\ref{phipipower} that 
    there
    exists some $L\gg 0$ such that $b^{(N)}_i\in\pp_n^L$ for all $i$. Since $\pp_n^L$ is an ideal, it
    follows that $x_i\in \pp_N^L$ for all $0\leq i\leq p-1$, which finishes the proof.
   \end{proof}

%++++++++++++++++++++++++++

   \subsection{Modules of differentials}     
  
     We make the following definitions:-
     \vs

     \noindent {\bf Definition.} Let 
     $\Omega^1_{\fSKdaggern\slash\ZZ}$ be the module of absolute differential forms of $\fSKdaggern$. 
     For $r\geq 1$, let $\Omega^r_{\fSKdaggern\slash\ZZ}=
     \bigwedge^r_{\fSKdaggern}\Omega^1_{\fSKdaggern\slash\ZZ}$,
     and define $\OmegarfSKn$ to be its $p$-adic completion
     \begin{equation}
      \OmegarfSKn=\varprojlim\Omega^r_{\fSKdaggern\slash\ZZ}\slash p^m
      \Omega^r_{\fSKdaggern\slash\ZZ}.
     \end{equation}
     
     \noindent {\bf Definition.} Let 
     \begin{equation}
      \Omega^1_{\fSKdaggern\slash\ZZ}(\log)=\Omega^1_{\fSKdaggern\slash\ZZ}
      \oplus\fSKdaggern\otimes_{\ZZ}(\fTKdaggern)^\times\slash\fN,
     \end{equation}
     where $\Omega^1_{\fSKdaggern\slash\ZZ}$ is the module of absolute differential forms and $\fN$ is the
     submodule of the direct sum generated by elements $(-da,a\otimes a)$ for $a\in\fSKdaggern\cap
     (\fTKdaggern)^\times$. In $\OmegafSKnlog$, denote the class $(0,1\otimes a)$ by $d\log a$. 
     For $r\geq 1$, let $\Omega^r_{\fSKdaggern\slash\ZZ}(\log)=
     \bigwedge^r_{\fSKdaggern}\Omega^1_{\fSKdaggern\slash\ZZ}(\log)$),
     and define $\OmegarfSKnlog$ to be its $p$-adic completion
     \begin{equation}
      \OmegarfSKnlog=\varprojlim\Omega^r_{\fSKdaggern\slash\ZZ}(\log)\slash p^m
      \Omega^r_{\fSKdaggern\slash\ZZ}(\log).
     \end{equation}
          
     \begin{lem}\label{structuredaggerdifferentials}
      We have the following isomorphisms of $\fSKdaggern$-modules:-
      \begin{align}
       \OmegatwofSKn&\cong \fSKdaggern . \bigwedge d\log(T_i)\wedge d\log (\pi_K+1),\\
       \OmegatwofSKnlog&\cong \fSKdaggern .\bigwedge d\log(T_i)\wedge d\log (\pi_K).
      \end{align}
     \end{lem}
     \begin{proof}
      Follows from explicit computation, using Lemma~\ref{formalderivative} and the observation 
      that $\pi_K+1$ is invertible in $\fSKdaggern$.
     \end{proof}
     
     Now $\fS_K^{\dagger,n+1}$ is finite flat and locally of complete intersection over
     $\phi(\fSKdaggern)$, 
     so as shown in~\cite{fukaya1}, we have a trace map of differential forms
     \begin{align}
      \Tr^{(n)}_\phi:&\Omega^n_{\fS_K^{\dagger,N+1}}\rTo \OmegatwofSKn.
     \end{align}
     and the usual trace map
     \begin{equation}
      \Tr_\phi: \fS_K^{\dagger,N+1}\rTo \fSKdaggern.
     \end{equation}
     \vs

     \begin{lem}\label{explicitdaggertrace}
      The trace map is characterized as follows:-
      \begin{equation}
       \Tr^{(n)}_\phi(c. \bigwedge d\log(T_i)\wedge d\log (\pi_K+1))
       =\frac{1}{p^2}\Tr_\phi (p.c\omega) \bigwedge d\log(T_i)\wedge d\log (\pi_K+1),
      \end{equation}
      where 
      \begin{equation}
       \omega=\frac{\phi(\pi_K)+1}{\pi_K+1}(\frac{d}{d\pi_K}\phi(\pi_K+1))^{-1}.
      \end{equation}
     \end{lem}
     \begin{proof}
      First note that since $\phi$ is a ring homomorphism, we have
      $\frac{d}{d\pi_K}\phi(\pi_K+1)=\frac{d}{d\pi_K}\phi(\pi_K)$. Now by
      Corollary~\ref{structurederivative}, there exists $I<0$ and $a\in(\AA_K^{\dagger,M})^\times$ for 
      $M\gg 0$ such that $$\frac{d}{d\pi_K}\phi(\pi_K)=p\pi_K^Ia,$$ so $p\omega\in\fS_K^{\dagger,M}$ for
      $M\gg 0$. The lemma is now an immediate consequence of Remark (iii) in Section 2.1
      in~\cite{fukaya1}.
     \end{proof}
     
     Using these formulae, we extend the trace map to $\OmegatwofSKnlog$:-
     \begin{equation}
      \Tr^{(n)}_\phi(\frac{c}{\pi_K}. \bigwedge d\log(T_i)\wedge  d\log (\pi_K+1))
      =\frac{1}{p^n}\frac{\Tr_\phi (p.c\tilde{\omega})}{\pi_K} \bigwedge d\log(T_i)\wedge 
      d\log (\pi_K+1),
     \end{equation}
     where
     \begin{equation}
      \tilde{\omega}=\frac{\phi(\pi_K)}{\pi_K}\frac{\phi(\pi_K)+1}{\pi_K+1}(\frac{d}{d\pi_K}\phi(\pi_K+1))^{-1}.
     \end{equation}
          
     \noindent {\bf Definition.} Let $\psi:\OmegatwofSKnlog\rTo\OmegatwofSKnlog$ be the map obtained by
     composing the natural map $\OmegatwofSKnlog\rTo\Omega^n_{\fS_K^{\dagger,N+1}}(\log)$ 
     with $\Tr^{(n)}_\phi$.

%++++++++++++++++++++++++++++++

  \section{The logarithmic derivative revisited}\label{Logderivativerevisited}
  
   \begin{lem}\label{injectionofdifferentials}
    The natural map $\iota: \OmegatwofSKnlog\rTo\OmegatwofRKlog$ is a closed embedding.
   \end{lem}
   \begin{proof}
    Recall that 
    \begin{align}
     \OmegatwofSKnlog&\cong\fSKdaggern d\log(\pi_K)\wedge d\log (T),\\
     \OmegatwofRKlog&\cong\fR_K^+ d\log(\pi_K)\wedge d\log (T)
    \end{align}
    The map $\iota$ is just the natural inclusion $\fSKdaggern\hookrightarrow\fR_K^+$, which is a closed
    embedding.
   \end{proof}
   
   As a corollary, we get the desired information about the logarithmic derivative:-
   
   \begin{thm}\label{logderivrevisited}
    The logarithmic derivative gives a canonical homomorphism
    \begin{equation}
     \Tau:\varprojlim\whK_n(\fK_i)\rTo(\OmegatwofSKnlog)^{\psi=1}.
    \end{equation}
   \end{thm}
   \begin{proof}
    It is certainly sufficient to show that we get a canonical homomorphism
    \begin{equation}
     \whK_n(\EfK)\rTo (\OmegatwofSKnlog)^{\psi=1}.
    \end{equation}
    Now if $\tilde{x}\in\whK_n(\EfK)$, then since $n\geq n(K)$ we can chose a lift $x$ of $\tilde{x}$ 
    such that $x$ can be written as a product of symbols $\{b_1,\dots,b_n\}$ such that 
    $b_i\in\fTKdaggern$ for all $i$. Then $d\log (x)\in\OmegatwofSKnlog$, and the sequence
    $(\psi^k(d\log(x)))$
    then converges to $y=\tau(\tilde{x})\in(\OmegatwofSKnlog)^{\psi=1}$. But $\psi$
    restricts to $\image(\iota)$, so since 
    $\OmegatwofSKnlog\rTo\OmegatwofRKlog$ is a closed embedding, it follows that 
    $y\in\OmegatwofSKnlog$, which finishes the proof.
   \end{proof}
   
   We have the following special case when $K=\Qp$:-
   
   \begin{cor}
    Suppose that $\fK=\fH$, i.e. $K=\Qp$. Then the logarithmic derivative gives a canonical homomorphism
    \begin{equation}
     \Tau:\varprojlim\whK_n(\fH_n)\rTo(\Omega^n_{O_{\fH}[[\pi]]}(\log))^{\psi=1}.
    \end{equation}
   \end{cor}
   
   \noindent {\bf Remark.} Note that 
   \begin{equation}
    \Omega^n_{O_{\fH}[[\pi]]}(\log)\cong O_{\fH}[[\pi]].\bigwedge d\log(T_i)\wedge d\log(\pi).
   \end{equation}

%++++++++++++++++++++++++++++++++++++++++++

  \section{Relation to the dual exponential map}\label{relationdualexp}
  
   \subsection{Statement of the theorem}
  
    Let $N\geq n(K)$.
    In this section we will use the logarithmic derivative to give a new descripton of Kato's dual
    exponential map $\lambda_m$ for $m\gg N$. 
    Recall that for all $m\geq 1$, the homomorphism $\lambda_m$ is defined as
    the composition
    \begin{align}
     \lambda_m: \varprojlim \whK_n(\fK_i)&\rTo^{(\delta^n_{\fK_i,p^i})_{i\geq 1} }
                                          \varprojlim H^n(\fK_i,(\ZZ\slash p^i\ZZ)(n))\\
                                        &\rTo^{\zeta_{p^i}^{\otimes(-1)}}
                                         \varprojlim H^n(\fK_i,(\ZZ\slash p^i\ZZ)(n-1))\\
                                        &\rTo^{\trace}\varprojlim H^n(\fK_m,(\ZZ\slash p^i\ZZ)(n-1))\\
                                        &=H^n(\fK_m,\ZZ_p(n-1))\\
                                        &\rTo H^n(\fK_m,\CC_p(n-1))\cong \Omega^{n-1}_{\fK_m},
    \end{align}
    where $\delta^n_{\fK_i,p^i}$ is the Galois symbol. Now let $\del_m$ be the homomorphism
    \begin{align}
     \del_m: \varprojlim \whK_n(\fK_i)&\rTo^{\Tau} (\OmegatwofSKnlog)^{\psi=1}\\
                                     &\hookrightarrow \OmegatwofSKnlog\cong \frac{\fSKdaggern}{\pi_K} .
                                     \bigwedge d\log(T_i)\wedge d\log (\pi_K+1)\\
                                     &\rTo \Omega^{n-1}_{\fK_m}.
    \end{align}
    Here the last arrow is defined by
    \begin{equation}
     f.\bigwedge d\log(T_i)\wedge d\log (\pi_K+1)\rTo \frac{1}{p^m}h^{(N)}_m(f).\bigwedge 
     d\log(h^{(N)}_m(T_i)) 
    \end{equation}
    where
    \begin{equation}
     h^{(N)}_m:\fSKdaggern\rightarrow O_{\fK_m}
    \end{equation}
    is determined by $h^{(N)}_m(T_i)=T_i^{\frac{1}{p^m}}$ and $h^{(N)}_m(\pi_K)=\pi_{K_m}$. 
    (Note that
    this is well-defined for $m\geq N$ since the elements of $\fSKdaggern$ can be interpreted as functions
    converging on the annulus $\{ p^{-\frac{1}{e(p-1)p^N}}\leq \mid X\mid <1\}\subset\CC_p$.) 
    We then have the following result:-

    \begin{thm}\label{dualexp}
     For all $m\geq N$, we have
     \begin{equation}
      \del_m=-\lambda_m.
     \end{equation}
    \end{thm} 

%+++++++++++++++++++++++++++

   \subsection{Proof}
   
    \subsubsection{Evaluation of the logarithmic derivative}
    
    \noindent {\bf Definition.} 
    For a complete discrete valuation field $L$ of mixed characteristic $(0,p)$, let
    \begin{equation}
     \Omega^1_{O_L\slash\ZZ}(\log)=(\Omega^1_{O_L\slash\ZZ}\oplus (O_L\otimes_{\ZZ}L^\times))\slash\fN,
    \end{equation}
    where $\fN$ is the $O_L$-submodule of the direct sum generated by elements of the form $(-da,a\otimes
    a)$for $a\in O_L\backslash\{ 0\}$. For $a\in L^\times$, denote the element $(0,1\otimes a)$ by $d\log
    (a)$. For $r\geq 1$, let 
    \begin{equation}
     \Omega^r_{O_L\slash\ZZ}(\log)=\bigwedge_{O_L}^r\Omega^1_{O_L\slash\ZZ}(\log),
    \end{equation}
    and define
    \begin{equation}
     \Omega^r_{O_L}(\log)=\varprojlim \Omega^r_{O_L\slash\ZZ}(\log)\slash p^i\Omega^r_{O_L\slash\ZZ}.
    \end{equation}
    
    Milnor $K$-theory then gives natural maps
    \begin{equation}
     d\log: K_r(L)\rTo\Omega^r_{O_L}(\log)
    \end{equation}
    which are defined by
    \begin{equation}
     \{a_1,\dots,a_r\}\rTo d\log(a_1)\wedge\dots\wedge d\log(a_n).
    \end{equation}
       
    For all $i\geq j\gg 0$, define the evaluation map
    \begin{equation}
     h^{(N)}_{i,j}:\fSKdaggern\rTo \fK_j(\zeta_{p^i})
    \end{equation}
    determined by $h^{(N)}_{i,j}(\pi_K)=\pi_{K_i}$ and $h^{(N)}_{i,j}(T_k)=(T_k^{\frac{1}{p^j}})$ for all $k$.
    Note that
    $h^{(N)}_{j,j}=h^{(N)}_{j}$, which was defined above.
    Let 
    \begin{align}
     H^{(N)}_{i,j}:& \whK_n(\fSKdaggern)\rTo \whK_n(\fK_j(\zeta_{p^i}))\\
     \fH_{i,j}^{(N)}:&\Omega^n_{\fSKdaggern}(\log)\rTo \Omega^n_{O_{\fK_i(\zeta_{p^i})}}(\log)
    \end{align}
    be the (surjective) maps induced by the evaluation map $h_{i,j}^{(N)}$. If $i=j$, denote them
    by $ H^{(N)}_i$ and $\fH_{i}^{(N)}$.
   
   \begin{prop}\label{evaluation}
    We have a commutative diagram
    \begin{diagram}
     \varprojlim\whK_n(\fK_i)   &  \rTo^{d\log}                   &  \varprojlim\Omega^n_{O_{\fK_i}}(\log)\\
     \dTo^\fT                       &  \ruTo_{(\fH^{(N)}_i)_{i\geq 1}}  & \\                
     (\OmegatwofSKnlog)^{\psi=1}&                                 & \\
    \end{diagram}
   \end{prop}
   \begin{proof}
    Recall that we have a continuous homomorphism $g:\varprojlim
    \whK_n(\fK_i)\rightarrow\whK_n(\EfK)$. Also, recall the following construction from the proof of
    Proposition~\ref{isoLCFT}:-
    For all $i\gg 0$, the evaluation map $h^{(N)}_i$ induces a homomorphism
    \begin{equation}
     H^{(N)}_i:\whK_n(\fTKdaggern)\rTo \whK_n(\fK_i).
    \end{equation}
    One can show that there exists $j_i^{(N)}\geq i-N-1$ such that 
    the image of $\ker(H^{(N)}_i)$ in $\whK_n(\EfK)$ is the
    subgroup $\fA^{(N)}_i$ generated by the symbols of the form $\{ a_1,\dots,a_n\}$ with 
    $a_1\in 1+\bpi_K^{j_i^{(N)}}\EfK^+$, and we have commutative diagram
    \begin{diagram}
     K_n(\EfK)\slash\fA^{(N)}_{i+1}&\lTo &\whK_n(\fTKdaggern)\slash\ker(H^{(N)}_{i+1}) 
     &\rTo^\cong_{H^{(N)}_{i+1}}&   \whK_n(\fK_{i+1})    \\
       \dTo                        &     &                       \dTo_N                      &            &    \dTo_{N_{O_{\fK_{m+1}}\slash O_{\fK_m}}}    \\
     K_n(\EfK)\slash\fA^{(N)}_i    &\lTo &\whK_n(\fTKdaggern)\slash\ker(H^{(N)}_i)     
     &\rTo^{\cong}_{H^{(N)}_i}  &   \whK_n(\fK_i)\\
    \end{diagram} 
    Now $j_i^{(N)}\rightarrow +\infty$ as $i\rightarrow+\infty$, so we have
    $\varprojlim K_n(\EfK)\slash\fA^{(N)}_i =\whK_n(\EfK)$, and
    the above diagram induces a map
    \begin{equation}
     G:\varprojlim \whK_n(\fTKdaggern)\slash\ker(H^{(N)}_i)\rTo\whK_n(\EfK).
    \end{equation}
    Now for all $i\gg 0$ we have a commutative diagram
    \begin{diagram}
     \fH_{i+1}^{(n)}:&\Omega^n_{\fSKdaggern}(\log)&\rTo&\Omega^n_{O_{\fK_{i+1}}}(\log)\\
                     & \dTo_{\Tr^{(n)}_\phi}      &          & \dTo_{\Tr}\\
     \fH_i^{(n)}:&\Omega^n_{\fSKdaggern}(\log)&\rTo&\Omega^n_{O_{\fK_i}}(\log)\\
    \end{diagram}
    so
    $\Tr^{(n)}_\phi(\ker(\fH_{i+1}^{(N)}))\subset\ker(\fH_i^{(N)})$ for all $i$, and we get 
    a commutative diagram
    \begin{diagram}
     G:&\varprojlim \whK_n(\fTKdaggern)\slash\ker(H^{(N)}_i)  &  \rTo^\cong & \varprojlim\whK_n(\fK_i)\\
       & \dTo_{d\log}                                        &         & \dTo_{d\log}\\
     (\fH_i^{(n)})_{i\gg 0}:&\varprojlim \Omega^n_{\fSKdaggern}(\log)\slash\ker(\fH_i^{(N)})&\rTo^\cong&
     \varprojlim\Omega^n_{O_{\fK_i}}(\log)\\
    \end{diagram}
    Now let $x\in\varprojlim\whK_n(\fTKdaggern)\slash\ker(H^{(N)}_i)$, and denote by $\bar{x}$ its image in
    $\whK_n(\EfK)$. Choose a sequence $(y^{(j)})$ of
    elements in $\tK_n(\fS_K)$ such that for all $i\gg 0$ 
    we have $\lambda(y^{(1)}\dots y^{(i)})\cong \bar{x}\mod\fA_i^{(N)}$.
    Note that we can choose $y^{(j)}$ such that it can
    be written as a product of symbols with entries in $\fTKdaggern\subset\fS_K$. 
    We will show that for all $i\gg 0$, $\psi^k(d\log (y^{(1)}\dots y^{(i+k)}))
    \rightarrow d\log(x)\mod \fH_i^{(N)}$ as
    $k\rightarrow +\infty$. Since $\fH_i^{(N)}$ is the evaluation map on the differentials, this will
    prove the proposition. 
    Now for each $j\gg 0$ there exist $y_1^{(j)},y_2^{(j)}\in\whK_n(\fS_K)$ 
    such that $y^{(1)}_1\dots y^{(i)}_1$ is a lift of 
    $x\mod  \ker(H^{(N)}_i)$ and 
    $y^{(1)}_2\dots y^{(i)}_2\in\ker(\lambda)$. By the definition of $y_1^{(j)}$, certainly 
    $\psi^kd\log(y_1^{(1)}\dots y_1^{(i+k)})=d\log(x)\mod\fH_i^{(N)}$. 
    But $N^k(y_2^{(1)}\dots y_2^{(i+k)})\rightarrow 1$ as $k\rightarrow +\infty$
    by Lemma~\ref{nilpotentaction}, so $\psi^k(d\log(y_2^{(1)}\dots y_2^{(i+k)}))\rightarrow 0$ 
    as $k\rightarrow +\infty$,
    which finishes the proof.
   \end{proof}
   
%+++++++++++++++++++++
   
   \subsubsection{The maps $\calf_{(\zeta_{p^i})_{i^\geq 1},m}$}

    Define 
    \begin{align}
     \lambda_{(\zeta_{p^i})_{i^\geq 1},m}: \varprojlim \whK_n(\fK_m(\zeta_{p^i}))
     &\rTo^{(\delta^n_{\fK_m(\zeta_{p^i}),p^i})_{i\geq 1} }
                                          \varprojlim H^n(\fK_m(\zeta_{p^i}),(\ZZ\slash p^i\ZZ)(n))\\
                                        &\rTo^{\zeta_{p^i}^{\otimes(-1)}}
                                         \varprojlim H^n(\fK_m(\zeta_{p^i}),(\ZZ\slash p^i\ZZ)(n-1))\\
                                        &\rTo^{\trace}\varprojlim H^n(\fK_m,(\ZZ\slash p^i\ZZ)(n-1))\\
                                        &=H^n(\fK_m,\ZZ_p(n-1))\\
                                        &\rTo H^n(\fK_m,\CC_p(n-1))\cong \Omega^{n-1}_{\fK_m},
    \end{align}
    In~\cite{kato3}, Kato defines a homomorphism
    \begin{equation}
     \calf_{(\zeta_{p^i})_{i^\geq 1},m}:\varprojlim \Omega^n_{O_{\fK_m(\zeta_{p^i}))}}(\log)
     \rTo \Omega^{n-1}_{\fK_m},
    \end{equation}
    and he shows the following result, which is cucial for our description of the dual exponential map:-
    
    \begin{thm}\label{katothm}
     For $m\geq 1$, let $\del_{(\zeta_{p^i})_{i^\geq 1},m}=\calf_{(\zeta_{p^i})_{i^\geq 1},m}\circ d\log$
     be the map
     \begin{align}
      \varprojlim \whK_n(\fK_m(\zeta_{p^i}))&\rTo^{d\log}\varprojlim 
      \Omega^n_{O_{\fK_m(\zeta_{p^i}))}}(\log)\\
      &\rTo^{\calf_{(\zeta_{p^i})_{i^\geq 1},m}}\Omega^{n-1}_{\fK_m}.
     \end{align}
     Then 
     \begin{equation}
      \lambda_{(\zeta_{p^i})_{i^\geq 1},m}=\del_{(\zeta_{p^i})_{i^\geq 1},m}.
     \end{equation}
    \end{thm}
    
    In the rest of this subsection, we will give a construction of the homomorphism
    $\calf_{(\zeta_{p^i})_{i^\geq 1},m}$, which we quote from Secton 5 in~\cite{fukaya2}:-
    
    \begin{lem}\label{katolem1}
     For all $i\geq m$, we have
     \begin{equation}
      p^{m-i}\Tr_{\fK_m(\zeta_{p^i})\slash\fK_m}(O_{\fK_m(\zeta_{p^i})})\subset O_{\fK_m}
     \end{equation}
     where $\Tr_{\fK_m(\zeta_{p^i})\slash\fK_m}:\fK_m(\zeta_{p^i})\rightarrow\fK_m$ is the usual trace
     map.
    \end{lem}
    \begin{proof}
     Lemma 6.1.6 in~\cite{kato3}.
    \end{proof}
    
    \begin{lem}\label{katolem2}
     Let $i\geq m$, and define
     \begin{align}
      \kappa_i:  O_{H_m(\zeta_{p^i})}\slash (p^i)\otimes_{O_\fK}\Omega^1_{O_{\fK}}&\rTo
      & \Omega^2_{O_{\fK_m(\zeta_{p^i})}}(\log) \\
        x\otimes y& \rTo & xy\wedge d\log(\zeta_{p^i})
     \end{align}
     Then $p^{e+m+2}$ kills $\ker(\kappa_i)$ and $p^{e+m+1}$ kills $\coker(\kappa_i)$.
    \end{lem}
    \begin{proof}
     Lemma 6.1.7 in~\cite{kato3}.
    \end{proof}
    
    Now let $c=c_1c_2c_3$ for $c_1=p^m$, $c_2=p^{e+m+1}$ and $c_3=p^{e+m+2}$, and define
    \begin{equation}
     f_{i,c}: \whK_n(\fK_m(\zeta_{p^i}))\rTo O_{\fK_m}\slash p^i\otimes_{O_{\fK}}\Omega^{n-1}_{O_{\fK}}
    \end{equation}
    as the composition
    \begin{align}
     \whK_n(\fK_m(\zeta_{p^i}))&\rTo^{d\log} \Omega^n_{\fK_m(\zeta_{p^i})}(\log)\\
                               &\rTo^{c_2\kappa_i^{-1}c_3}
                               O_{\fK_m(\zeta_{p^i})}\slash p^i\otimes_{O_{\fK}}\Omega^{n-1}_{O_{\fK}}\\
                               &\rTo^{c_1p^{-i}\Tr_{\fK_m(\zeta_{p^i})\slash\fK_m}\otimes\id}
                               O_{\fK_m}\slash p^i\otimes_{O_{\fK}}\Omega^{n-1}_{O_{\fK}}
   \end{align}
   Here, the map $c_2\kappa_i^{-1}c_3$ is defined as follows:- For
   $x\in\Omega^n_{\fK_m(\zeta_{p^i})}(\log)$, let $y$ be an element of 
   $O_{\fK_m(\zeta_{p^i})}\slash p^i\otimes_{O_{\fK}}\Omega^{n-1}_{O_{\fK}}$ such that $c_3x=\kappa_i(y)$.
   Then $c_2\kappa_i^{-1}c_3(x)=c_2(y)$.
   \vs
   
   We can now define the map $\calf_{(\zeta_{p^i})_{i^\geq 1},m}$:-
   \vs
   
   \noindent {\bf Definition.} Let 
   \begin{equation}
    \calf_{(\zeta_{p^i})_{i^\geq 1},m}=c^{-1}\varprojlim_i f_{i,c}.
   \end{equation}

%+++++++++++++++++++

   \subsubsection{Lots of commutative diagrams}

    In view of Theorem~\ref{katothm}, 
    the proof of Theorem~\ref{dualexp} reduces to showing the following result:-
    
    \begin{prop}\label{relationdels}
     The composition
     \begin{align}\label{defnmap}
      \varprojlim\whK_n(\fK_i)&\rTo^{(N_{\fK_i\slash \fK_m(\zeta_{p^i})})_i} 
      \varprojlim \whK_n(\fK_m(\zeta_{p^i}))\\
      &\rTo^{\calf_{(\zeta_{p^i})_{i^\geq 1},m}}\Omega^{n-1}_{\fK_m}
     \end{align}
     coincides with $\del_m$.
    \end{prop}
    
    Using the commutative diagram
    \begin{diagram}
     \whK_n(\fK_i)                        &  \rTo^{d\log}  & \Omega^n_{O_{\fK_i}}(\log) \\
     \dTo_{N_{\fK_i\slash \fK_m(\zeta_{p^i})}}&         & \dTo_{\Tr_{\fK_i\slash \fK_m(\zeta_{p^i})}}\\
     \whK_n(\fK_m(\zeta_{p^i}))    &  \rTo^{d\log}  & \Omega^n_{O_{\fK_m(\zeta_{p^i})}}(\log) \\
    \end{diagram}
    we see that
    to prove Proposition~\ref{relationdels} it is sufficient to show the following result:-
    
    \begin{prop}
     The map 
     \begin{align}
     \varprojlim\whK_n(\fK_i)&\rTo^{d\log}\varprojlim\Omega^n_{O_{\fK_i}}(\log)\\
                             &\rTo^{\Tr_{\fK_i\slash \fK_m(\zeta_{p^i})}}
                                 \varprojlim\Omega^n_{O_{\fK_m(\zeta_{p^i})}}(\log)\\
                             &\rTo^{\calf_{(\zeta_{p^i})_{i^\geq 1},m}}\Omega^{n-1}_{\fK_m}
    \end{align}
    agrees with $\del_m$.
   \end{prop}
   \begin{proof}
    Define homomorphisms $\sigma,\tau: \fSKdaggern\rightarrow \fS_K^{\dagger,n+1}$, which are determined by
    \begin{align}
     \sigma:& \pi_K\rTo\phi(\pi_K),\\
     \tau:& T_i \rTo T_i^p\text{    for all $i$}.
    \end{align}
    Note that $\phi=\sigma\circ\tau$. As before, these homomorphisms then induce trace maps
    \begin{equation}
     \Tr^{(n)}_\sigma,\Tr^{(n)}_\tau: \Omega^n_{\fSKdaggern}(\log)\rTo\Omega^n_{\fSKdaggern}(\log).
    \end{equation}
    Explicit computation shows that 
    \begin{align*}
     \Tr^{(n)}_\sigma(\frac{a}{\pi_K}\bigwedge d\log(T_i)\wedge d\log(\pi_K+1))&
      =\frac{1}{p}\frac{\Tr_\sigma(p.a\omega)}{\pi_K} \bigwedge d\log(T_i)\wedge 
      d\log (\pi_K+1)\\
     \Tr^{(n)}_\tau(\frac{b}{\pi_K}\bigwedge d\log(T_i)\wedge d\log(\pi_K+1))&
      =\frac{1}{p^{n-1}}\frac{\Tr_\tau(p.b\omega)}{\pi_K} \bigwedge d\log(T_i)\wedge 
      d\log (\pi_K+1)
    \end{align*}
    where
    \begin{equation}
     \omega=\frac{\phi(\pi_K)}{\pi_K}\frac{\phi(\pi_K)+1}{\pi_K+1}(\frac{d}{d\pi_K}\phi(\pi_K+1))^{-1},
    \end{equation}
    and we have a commutative diagram
    \begin{diagram}
     \Omega^n_{\fSKdaggern}(\log)  &\rTo^{(\Tr^{(n)}_\sigma)^{i-m}} & \Omega^n_{\fSKdaggern}(\log)\\
            \dTo_{\fH^{(N)}_{i,m}}  &     & \dTo_{\fH^{(N)}_m}\\                   \\
     \Omega^n_{O_{\fK_m(\zeta_{p^i})}}(\log) &\rTo^{\Tr^{(n)}_{\fK_i\slash\fK_m(\zeta_{p^i})}} 
     & \Omega^n_{O_{\fK_i}}(\log)\\
    \end{diagram}
    Also, $\tau$ induces a trace map 
    $$\Tr_\tau:\frac{\fSKdaggern}{\pi_K}\rTo\frac{\fSKdaggern}{\pi_K},$$
    and we have a commutative diagram
    \begin{diagram}
     \frac{\fSKdaggern}{\pi_K} & \rTo^{(\Tr_\tau)^{i-m}} & \frac{\fSKdaggern}{\pi_K}\\
      \dTo_{h^{(N)}_{i,m}}       &                         & \dTo_{h^{(N)}_i}\\
     \fK_m(\zeta_{p^i})       & \rTo^{\Tr_{\fK_m(\zeta_{p^i})\slash\fK_m}}& \fK_m\\
    \end{diagram}
    Here, $\Tr_\tau(\frac{c}{\pi_K})=\frac{\Tr_\tau(\omega.c)}{\pi_K}$, where
    $\Tr_\tau:\fSKdaggern\rightarrow\fSKdaggern$ is the usual trace map. 
    
    Recalling the definition of $\calf_{(\zeta_{p^i})_{i^\geq 1},m}$  and combining the above 
    commutative diagrams with Proposition~\ref{evaluation} shows that the map 
    given by equations~\eqref{defnmap} and (144) is also given by the composition
    \begin{align}
     \varprojlim\whK_n(\fK_i)&\rTo &(\Omega^n_{\fSKdaggern}(\log))^{\psi=1}& &\\
                             &\hookrightarrow &\Omega^n_{\fSKdaggern}(\log)&
                              \rTo^{(\Tr^{(n)}_\sigma)^{i-m}} &\Omega^n_{\fSKdaggern}(\log)\\
                             & & &\rTo^{\ppp}&\frac{\fSKdaggern}{\pi_K}\bigwedge d\log(T_i)\\
                             & & &\rTo^{(\Tr_\tau)^{i-m}}&\frac{\fSKdaggern}{\pi_K}\bigwedge d\log(T_i)\\
                             & & &\rTo& \Omega^{n-1}_{\fK_m}
    \end{align}
    where the map $\ppp$
    is given by $$\ppp:f\bigwedge d\log(T_i)\wedge d\log(\pi_K+1)\rTo f\bigwedge d\log(T_i)$$ and the last arrow
    is defined as 
    \begin{equation}
     f\bigwedge d\log(T_i)\rTo \frac{1}{p^m}h^{(N)}_m(f)\bigwedge d\log(h^{(N)}_m(T_i)).
    \end{equation}
    But $\Tr^{(n)}_\tau\circ\Tr^{(n)}_\sigma=\Tr^{(n)}_\phi$ and
    $\Tr_\tau(\ppp(x))=\ppp(\Tr^{(n)}_\tau(x))$ for all $x\in\Omega^n_{\fSKdaggern}(\log)$, which 
    shows that this composition also gives the map $\del_m$.
   \end{proof}

 %++++++++++++++++++++++++  

\end{document}